\documentclass[preprint,12pt,compress]{elsarticle}

\usepackage{amssymb}
\usepackage{amsmath}
\usepackage[a4paper, body={15.6cm,23.5cm}]{geometry}
\usepackage{bbm}
\begin{document}
\journal{ArXiv}
\parindent 15pt
\parskip 4pt


 \newcommand{\eps}{\varepsilon}
 \newcommand{\lam}{\lambda}
 \newcommand{\To}{\rightarrow}
 \newcommand{\as}{{\rm d}P\times {\rm d}t-a.e.}
 \newcommand{\ps}{{\rm d}P-a.s.}
 \newcommand{\jf}{\int_t^T}
 \newcommand{\tim}{\times}

 \newcommand{\F}{\mathcal{F}}
 \newcommand{\E}{\mathbb{E}}
 \newcommand{\N}{\mathbb{N}}
 \newcommand{\s}{\mathcal{S}}
 \newcommand{\M}{{\rm M}}
 \newcommand{\T}{[0,T]}
 \newcommand{\Lp}{L^p(\Omega,\F_T, P)}

 \newcommand{\R}{{\mathbb R}}
 \newcommand{\Q}{{\mathbb Q}}
 \newcommand{\RE}{\forall}

 \newcommand{\US}{L^1([0,T];\R_+)}
 \newcommand{\VS}{L^2([0,T];\R_+)}
 \newcommand{\vS}{L^{2\over 2-\alpha}([0,T];\R_+)}
 \newcommand{\LS}{\mathcal{L}(\R;\R^+)}
 \newcommand{\FS}{L^p(\Omega;L^1(\T;\R_+))}

\newcommand{\gET}[1]{\underline{{\mathcal {E}}_{t,T}^g}[#1]}
\newcommand {\Lim}{\lim\limits_{n\rightarrow \infty}}
\newcommand {\Limk}{\lim\limits_{k\rightarrow \infty}}
\newcommand {\Limm}{\lim\limits_{m\rightarrow \infty}}
\newcommand {\llim}{\lim\limits_{\overline{n\rightarrow \infty}}}
\newcommand {\slim}{\overline{\lim\limits_{n\rightarrow \infty}}}
\newcommand {\Dis}{\displaystyle}

\begin{frontmatter}
\title {Bounded solutions, $L^p\ (p>1)$ solutions and $L^1$ solutions for one-dimensional BSDEs under general assumptions\tnoteref{fund}}
\tnotetext[fund]{Supported by the National Natural Science Foundation of China (No. 11371362), the China Postdoctoral Science Foundation (No. 2013M530173), the Qing Lan Project and the Fundamental Research Funds for the Central Universities (No. 2013RC20)\vspace{0.2cm}.}

\author{ShengJun FAN}
\ead{f$\_$s$\_$j@126.com}

\address{College of Sciences, China University of Mining and
Technology, Jiangsu 221116, P.R. China\vspace{-0.8cm}}

\begin{abstract}
This paper aims at solving one-dimensional backward stochastic differential equations (BSDEs) under weaker assumptions. We establish general existence, uniqueness, and comparison results for bounded solutions, $L^p\ (p>1)$ solutions and $L^1$ solutions of the BSDEs. The time horizon is allowed to be finite or infinite, and the generator $g$ is allowed to have a general growth in $y$ and a quadratic growth in $z$. As compensation, the generator $g$ needs to satisfy a kind of one-sided linear or super-linear growth condition in $y$, instead of the monotonicity condition in $y$ as is usually done. Many of our results improve virtually some known results, even though for the case of the finite time horizon and the case of the $L^2$ solution. \vspace{0.1cm}
\end{abstract}

\begin{keyword}
Backward stochastic differential equation\sep Existence and uniqueness \sep Comparison theorem\sep $L^p$ solution\sep Bounded solution\sep Infinite time interval\vspace{0.2cm}

\MSC[2010] 60H10
\end{keyword}
\end{frontmatter}\vspace{-0.3cm}


\section{Preliminaries and introduction}

Throughout this paper, let us fix a positive integer $d$ and a positive time horizon $T$ which can be finite or infinite. Moreover, let $\R_+:=[0,+\infty)$, $\R^+:=(0,+\infty)$ and let $\US$ and $\VS$, respectively, represent the set of nonnegative functions $u(\cdot)$ from $\T$ to $\R_+$ such that
$$\int_0^T u(t)\ {\rm d}t<+\infty\ \ {\rm and}\ \ \int_0^T u^2(t)\ {\rm d}t<+\infty.\vspace{-0.1cm}$$

Let $(\Omega,\F,P)$ be a complete probability space carrying a standard $d$-dimensional Brownian motion $(B_t)_{t\geq 0}$. We consider $(\F_t)_{t\geq 0}$ the natural filtration of $(B_t)_{t\geq 0}$ augmented by the $P$-null sets of $\F$ and we assume $\F_T=\F$. $\mathcal P$ denotes, as usual, the $\sigma$-algebra of predictable subsets of $\Omega\times \T$. For each predictable subset $A$ of $\Omega\times \T$, let $\mathbbm{1}_A$ equal to $1$ when $(t,\omega)\in A$, and $0$ otherwise. For every positive integer $n$, we use $| \cdot |$ to denote the norm of Euclidean space $\R^n$. For each $p>0$, $\Lp$ represents the set of $\F_T$-measurable random variable $\xi$ such that $\E[|\xi|^p]<+\infty$, and ${\s}^p$ denotes the set of real-valued, $(\F_t)$-adapted and continuous processes $(Y_t)_{t\in\T}$
such that
$$\|Y\|_{{\s}^p}:=\left( \E[\sup_{t\in\T} |Y_t|^p] \right)
^{1\wedge 1/p}<+\infty. $$
If $p\geq 1$, $\|\cdot\|_{{\s}^p}$ is a norm on ${\s}^p$ and if $p\in (0,1)$, $(X,X')\longmapsto \|X-X'\|_{{\s}^p}$ defines a distance on
${\s}^p$. Under this metric, ${\s}^p$ is complete. Moreover, for each $p>0$, let ${\rm M}^p$ denote the set of (equivalent classes of)
$(\F_t)$-progressively measurable, ${\R}^d$-valued processes $(Z_t)_{t\in\T}$ such that
$$\|Z\|_{{\rm M}^p}:=\left\{ \E\left[\left(\int_0^T |Z_t|^2\
{\rm d}t\right)^{p/2}\right] \right\}^{1\wedge 1/p}<+\infty.
$$
For $p\geq 1$, ${\rm M}^p$ is a Banach space endowed
with this norm and for $p\in (0,1)$, ${\rm M}^p$ is a complete metric space with the resulting distance. Finally, for each $p>1$, we denote by $\FS$ the set of $(\F_t)$-progressively measurable, nonnegative and real-valued process $(f_t)_{t\in\T}$ such that
$$
\E\left[\left(\int_0^T |f_t|\ {\rm d}t\right)^p\right]<+\infty.
$$
And, $L^1(\Omega;L^1(\T;\R_+))$ will be denoted simply by $L^1(\Omega\times\T)$.
 
We set $\s=\cup_{p>1} {\s}^p$ and denote by $\s^\infty$ the set of predictable bounded processes, and by $L^\infty(\Omega,\F_T,P)$ the set of $\F_T$-measurable bounded random variables. Let us recall that a continuous process $(Y_t)_{t\in\T}$ belongs to the class (D) if the family $\{ Y_\tau:\tau\in\Sigma_T\}$ is uniformly integrable, where and hereafter $\Sigma_T$ stands for the set of
all $(\F_t)$-stopping times $\tau$ such that $\tau\leq T$.

In this paper, we consider the following one-dimensional backward stochastic differential equation (BSDE in short for the remaining of this paper):
\begin{equation}
y_t=\xi+\int_t^Tg(s,y_s,z_s){\rm d}s-\int_t^Tz_s\cdot {\rm d}B_s,\ \
    t\in\T,
\end{equation}
where the random variable $\xi$ is $\F_T$-measurable, called the terminal condition of BSDE(1), the random function
$$g(\omega,t,y,z):\Omega\tim \T\tim {\R}\tim {\R}^d \longmapsto \R$$
is ${\mathcal P}\otimes {\mathcal B}(\R)\otimes {\mathcal B}(\R^d)$ measurable, called the generator of BSDE(1). We will sometimes use the notation BSDE$(\xi,g)$ to say that we consider the BSDE whose generator is $g$ and whose terminal condition is $\xi$.

For convenience of the following discussion, we introduce the following definitions concerning the solutions of BSDE(1).\vspace{0.1cm}

{\bf Definition 1.1}\ A solution of BSDE(1) is a pair of $(\F_t)$-progressively measurable processes $(y_\cdot,z_\cdot)$ with values in ${\R}\times {\R}^d$ such that $\ps$, $t\mapsto y_t$ is continuous, $t\mapsto z_t$ belongs to ${\rm L}^2(0,T)$, $t\mapsto g(t,y_t,z_t)$ belongs to ${\rm L}^1(0,T)$, and $\ps$, (1) holds true for each $t\in\T$.\vspace{0.1cm}

{\bf Definition 1.2}\ Assume that $(y_\cdot,z_\cdot)$ is a solution of BSDE(1). If $(y_\cdot,z_\cdot)\in \s^\infty\times {\rm M}^2$, then it will be called a bounded solution; if $(y_\cdot,z_\cdot)\in\s^p\times {\rm M}^p$ for some $p>1$, then an $L^p$ solution; if $(y_\cdot,z_\cdot)\in \s^\beta\tim\M^\beta$ for any $\beta\in (0,1)$ and $y_\cdot$ belongs to the class (D), then an $L^1$ solution.\vspace{0.1cm}

{\bf Definition 1.3}\ We say that $(y_\cdot,z_\cdot)$ is a maximal bounded (resp. $L^p\ (p>1)$ and $L^1$) solution of BSDE(1) if it is a bounded (resp. $L^p$ and $L^1$) solution and $y_\cdot\geq y'_\cdot$ for any bounded (resp. $L^p$ and $L^1$) solution $(y'_\cdot,z'_\cdot)$. Similarly, if $y_\cdot\leq y'_\cdot$, then $(y_\cdot,z_\cdot)$ will be called a minimal bounded (resp. $L^p$ and $L^1$) solution.\vspace{0.2cm}

Since the first existence and uniqueness result for nonlinear multidimensional BSDEs with square integrable parameters was introduced by \citet{Par90} under the Lipschitz assumption of $g$, BSDEs have been extensively studied, and many applications have been found in mathematical finance, stochastic control, partial differential equations and so on (see \cite{Bah10,BC08,Delb11,Delb10,El97,Hu05,Kob00,Mor09,Par99, Peng91,Peng97,Tang94,Tang98} for details).

From the beginning, many authors attempted to improve the result of \cite{Par90} by weakening the Lipschitz hypothesis on $g$, see \cite{Bah02,Bah10,BC08,Bri03,
Bri06,Bri08,Bri07,Delb11,Delb10,El97,Fan10,Fan12,FJ12,
FJD10,FJT11,FL10,Ham03,Jia08,Jia10,Kob00,
Lep97,Lep98,Lep02,Mor09,Par99,Peng91,Tian13,Xiao12}, or the $L^2$ integrability assumptions on $\xi$, see \cite{Bri03,Bri06,Chen10,El97,Fan12,FJ12,FL10,Peng97,
Tian13,Xiao12}, or relaxing the finite time horizon $T$ to a stopping time or infinity, see \cite{Chen00,FJ12,
FJT11,Mor09,Par99}. From these results it is not difficult to see that the case of one-dimensional BSDEs is easier to handle due to the presence of the comparison theorem of solutions, see \cite{Bri06,Bri08,Cao99,Chen10,Chen00,El97,
Fan10,Fan12,FJT11,FL10,Jia08,Jia10,Kob00,Lep97,Lep98,
Lep02,Peng97,Tian13,Xiao12}.

This paper aims at solving one-dimensional BSDEs under weaker assumptions. We establish general existence, uniqueness, and comparison results for bounded solutions, $L^p\ (p>1)$ solutions and $L^1$ solutions of the BSDEs. The time horizon $T$ is allowed to be finite or infinite, and the generator $g$ is allowed to have a general growth in $y$ and a quadratic growth in $z$. As compensation, the generator $g$ needs to satisfy a kind of one-sided linear or super-linear growth condition in $y$, instead of the monotonicity condition in $y$ as is usually done. Many of our results improve virtually some known results, even though for the case of the finite time horizon and the case of the $L^2$ solution.

More specifically, with respect to the existence for bounded solutions, $L^p\ (p>1)$ solutions and $L^1$ solutions of BSDEs, we would like to list respectively several existing results and our results as follows in order to compare with each other. Here, we always assume that the generator $g$ is continuous in $(y,z)$, and $T$ is a finite real number.

First, when $\xi$ is bounded there exists a maximal (resp. minimal) bounded solution of BSDE$(\xi,g)$ under one of the following two groups of conditions:\vspace{-0.2cm}
\begin{itemize}
\item $g$ has a super-linear growth in $y$ and a quadratic growth in $z$, i.e., there exists a constant $C>0$ and a continuous function $l:\R\mapsto\R^+$ such that
    $$|g(\omega,t,y,z)|\leq l(y)+C|z|^2,$$
    where $\int_0^{+\infty} {{\rm d}x\over l(x)}=\int_{-\infty}^0 {{\rm d}x\over l(x)}=+\infty$\vspace{0.2cm}\\ (see \citet{Lep98} and \citet{Kob00}).
\item There exist two constants $\alpha>0$ and $\beta>0$ together with a convex and ${\mathcal C}^1$ function $\rho:\R_+\mapsto \R_+$ with $\rho(0)=0$ and a continuous function $\varphi:\R_+\mapsto \R_+$ with $\varphi(0)=0$ such \vspace{0.2cm}that
    $$\begin{array}{l}
    (i)\ \ \ (g(\omega,t,y,z)-g(\omega,t,0,z))\ {\rm sgn}(y)\leq \rho(|y|),\vspace{0.1cm}\\
    (ii) \ \ |g(\omega,t,y,z)|\leq \alpha+\varphi(|y|)+\beta |z|^2,\vspace{0.1cm}
    \end{array}$$
    where $\int_0^{+\infty} {{\rm d}s\over \rho(s)+\alpha}=+\infty$\vspace{0.2cm}\\
     (see \citet{Bri08} and \citet{Bri07}).\vspace{-0.1cm}
\end{itemize}
It is not hard to verify that neither of the above two groups of conditions is satisfied for the following two generators:
$$g_1(\omega,t,y,z):=|z|^2{\rm e}^y+y\cos y\ \ \ {\rm and}\ \ \  g_2(\omega,t,y,z):=-y^3+|z|^{3\over 2}\sin y.$$
However, they both satisfy the following condition:\vspace{-0.2cm}
\begin{itemize}
\item $g$ has a one-sided super-linear growth in $y$ and a quadratic growth in $z$, i.e., there exist two continuous functions $l:\R\mapsto\R^+$ and $\varphi:\R\mapsto\R_+$ such that
    $$\begin{array}{l}
    (i)\ \ \ g(\omega,t,y,z)\ {\rm sgn}(y)\leq l(y)+\varphi(y)|z|^2, \vspace{0.1cm}\\
    (ii) \ \ |g(\omega,t,y,z)|\leq \varphi(y)(1+|z|^2),
    \end{array}
    $$
    where $\int_0^{+\infty} {{\rm d}x\over l(x)}=\int_{-\infty}^0 {{\rm d}x\over l(x)}=+\infty$.
\end{itemize}
Then, by Theorem 3.1 of this paper we know that when $\xi$ is bounded there exists a maximal (resp. minimal) bounded solution of BSDE$(\xi,g_i)$ for $i=1,2$. In addition, Theorem 3.1 also considers the case of $T=+\infty$.\vspace{0.2cm}

Second, when $\xi\in \Lp$ for some $p>1$ there exists a maximal (resp. minimal) $L^p$ solution of BSDE$(\xi,g)$ under one of the following two groups of conditions:\vspace{-0.2cm}

\begin{itemize}
\item $g$ has a linear growth in $(y,z)$, i.e., there exists a constant $C>0$ such that
    $$|g(\omega,t,y,z)|\leq C(1+|y|+|z|)$$
    (see \citet{Lep97} and \citet{Chen10}).
\item There exist two constants $\mu\in \R$ and $A>0$ together with a continuous adapted process $g_t\in \FS$ and a continuous function $\varphi:\R_+\mapsto \R_+$ with $\varphi(0)=0$ such that
    $$\begin{array}{l}
    (i)\ \ \ (g(\omega,t,y_1,z)-g(\omega,t,y_2,z))\ {\rm sgn}(y_1-y_2)\leq \mu|y_1-y_2|,\vspace{0.1cm}\\
    (ii) \ \ |g(\omega,t,y,z)|\leq g_t(\omega)+\varphi(|y|)+A|z|
    \end{array}$$
    (see \citet{Bri07}).\vspace{-0.2cm}
\end{itemize}
It is not hard to verify that neither of the above two groups of conditions is satisfied for the following two generators:
$$g_1(\omega,t,y,z):=|z|^2(1-{\rm e}^y)+|z|\sin|z|\ \ \ {\rm and}\ \ \  g_2(\omega,t,y,z):=-y^5+\cos(y|z|).$$
However, they both satisfy the following condition:\vspace{-0.2cm}
\begin{itemize}
\item $g$ has a one-sided linear growth in $(y,z)$ and a quadratic growth in $z$, i.e., there exists a constant $C>0$ and a continuous function $\varphi:\R\mapsto \R_+$ such that
    $$\begin{array}{l}
    (i)\ \ \ g(\omega,t,y,z)\ {\rm sgn}(y)\leq C(1+|y|+|z|), \vspace{0.1cm}\\
    (ii) \ \ |g(\omega,t,y,z)|\leq \varphi(y)(1+|z|^2).
    \end{array}\vspace{-0.2cm}
    $$
\end{itemize}
Then, by Theorem 5.1 of this paper we know that when $\xi\in \Lp$ for some $p>1$ there exists an $L^p$ solution of BSDE$(\xi,g_i)$ for $i=1,2$. In addition, Theorem 5.1 also considers the case of $T=+\infty$, and Theorem 5.2 further investigates the existence of a maximal (resp. minimal) $L^p$ solution.\vspace{0.2cm}

Third, when $\xi\in L^1(\Omega,\F_T,P)$ there exists an $L^1$ solution of BSDE$(\xi,g)$ under one of the following two groups of conditions:\vspace{-0.2cm}
\begin{itemize}
\item $g$ has a linear growth in $y$ and a sub-linear growth in $z$, i.e., there exist two constants $C>0$ and $\alpha\in (0,1)$ such that
    $$|g(\omega,t,y,z)|\leq C(1+|y|+|z|^\alpha)$$
    (see the first version of \citet{Bri06}).
\item There exist constants $\mu\in \R$, $\lambda\geq 0$, $\delta\geq 0$ and $\alpha\in (0,1)$ together with an $(\F_t)$-progressively measurable nonnegative process $g_t\in L^1(\Omega\times \T)$ such that
    $$\begin{array}{l}
    (i)\ \ \ (g(\omega,t,y_1,z)-g(\omega,t,y_2,z))\ {\rm sgn}(y_1-y_2)\leq \mu|y_1-y_2|,\vspace{0.1cm}\\
    (ii) \ \ |g(\omega,t,y,z_1)-g(\omega,t,y,z_2)|\leq \lambda |z_1-z_2|,\vspace{0.1cm}\\
    (iii) \ {\rm for\ each\ } r\geq 0, \psi_r(t):=\sup\limits_{|y|\leq r}|g(\omega,t,y,0)|\in L^1(\Omega\times\T),\vspace{0.1cm}\\
    (iv) \ \ |g(\omega,t,y,z)-g(\omega,t,y,0)|\leq \delta (g_t(\omega)+|y|+|z|)^\alpha
    \end{array}$$
    (see \citet{Bri03} in the multidimensional \vspace{-0.25cm}case).
\end{itemize}
It is not hard to verify that neither of the above two groups of conditions is satisfied for the following two generators:
$$g_1(\omega,t,y,z):=-|z|^2y^3+\sqrt[3]{|z|}\ \ \ {\rm and}\ \ \  g_2(\omega,t,y,z):={\rm e}^{-y}\sqrt{|z|}+\sqrt{1+|y|+|z|}.$$
However, they both satisfy the following condition:\vspace{-0.2cm}
\begin{itemize}
\item $g$ has a one-sided linear growth in $y$, a one-sided sub-linear growth in $z$ and a quadratic growth in $z$, i.e., there exist two constants $C>0$, $\alpha\in (0,1)$ and a continuous function $\varphi:\R\mapsto \R_+$ such that\vspace{-0.1cm}
    $$\begin{array}{l}
    (i)\ \ \ g(\omega,t,y,z)\ {\rm sgn}(y)\leq C(1+|y|+|z|^\alpha), \vspace{0.1cm}\\
    (ii) \ \ |g(\omega,t,y,z)|\leq \varphi(y)(1+|z|^2).
    \end{array}\vspace{-0.3cm}
    $$
\end{itemize}
Then, by Theorem 6.1 of this paper we know that when $\xi\in L^1(\Omega,\F_T,P)$ there exists an $L^1$ solution of BSDE$(\xi,g_i)$ for $i=1,2$. In addition, Theorem 6.1 also considers the case of $T=+\infty$, and Theorem 6.2 further investigates the existence of a maximal (resp. minimal) $L^1$ solution.\vspace{0.2cm}

In the sequel, with respect to the uniqueness and comparison results for bounded solutions, $L^p\ (p>1)$ solutions and $L^1$ solutions of BSDEs, we also list respectively several existing results and our results as follows.

\citet{Bri08} established a comparison theorem for bounded solutions of BSDEs with finite time horizon when one of the generators is Lipschitz in $y$ and concave or convex in $z$; \citet{Mor09} obtained
a comparison theorem for bounded solutions of BSDEs with infinite time horizon when one of the generators satisfies a monotonicity condition in $y$ and a local Lipschitz condition in $z$. Under the conditions that one of the generators only satisfies a one-sided Osgood condition in $y$ (see assumption (2A1) and Remark 2.1 in Section 2 for details), and a local Lipschitz condition (resp. a concavity or convexity condition) in $z$, Theorem 4.1 and Theorem 4.2 of this paper respectively prove a comparison theorem for bounded solutions of BSDEs with finite or infinite time horizon, extending two results mentioned above.

\citet{FJT11} established a comparison theorem for $L^2$ solutions of BSDEs with finite or infinite time horizon when one of the generators satisfies a one-sided Osgood condition in $y$ and a uniform continuity condition in $z$, which generalizes four classical comparison theorems obtained respectively in \citet{El97}, \citet{Cao99}, \citet{Chen00} and \citet{Bri06}. Theorem 2.1 of this paper further extends this result to the case of the $L^p\ (p>1)$ solution. More importantly, in this paper we eliminate the concavity condition with respect to the function $\rho(\cdot)$ in the one-sided Osgood condition.

To our knowledge, \citet{Bri06} first put forward and prove a comparison theorem for $L^1$ solutions of BSDEs with finite time horizon when one of the generators satisfies a monotonicity condition in $y$ and a Lipschitz condition together with a sub-linear growth condition in $z$. Recently, \citet{FL10}, \citet{Xiao12}, \citet{Fan12} and \citet{Tian13} further establish the comparison results for $L^1$ solutions of BSDEs with finite time horizon under the conditions that one of the generators satisfies a monotonicity condition or a Osgood condition in $y$ and a quasi-H\"{o}lder continuity condition in $z$. Theorem 2.4 of this paper unifies these results to the case of BSDEs with finite or infinite time horizon when one of the generators satisfies a one-sided Osgood condition in $y$ and a uniform continuity condition together with a sub-linear growth condition in $z$. For example, the following generator $g$ does not satisfy their conditions but satisfies our conditions:
$$
g(\omega,t,y,z)=\left\{
\begin{array}{rcl}
\sqrt{|z|}&,& 0\leq |z|\leq 1;\\
(n-1)^4 &,& (n-1)^4<|z|\leq n^4-2n+1;\\
|z|+n^2-n^4 &,& n^4-2n+1<|z|\leq n^4,
\end{array}\right.
\ \ \ \ n=2,3,4,\cdots.
$$
More specifically, this generator $g$ satisfies the uniform continuity condition as well as the sub-linear growth condition in $z$, but it does not satisfy the quasi-H\"{o}lder continuity condition in $z$.

Finally, we would like to mention that several new comparison results for the maximal and minimal solutions of BSDEs are put forward and proved (see, for example, Theorems 2.2-2.3, 3.2-3.3, 5.3-5.4 and 6.3-6.4), and they also play an important role in the proof of our existence results. Furthermore, by products, three general existence and uniqueness results for bounded solutions, $L^p\ (p>1)$ solutions and $L^1$ solutions of BSDEs are also obtained respectively (see Theorems 4.3, 5.5 and 6.5). In addition, we also point out that our results are all obtained due to the application of new ideas and new techniques or the development of those existing ideas and methods.

The remaining of this paper is organized as follows. In Section 2, we establish several comparison theorems for  the (maximal and minimal) $L^p\ (p> 1)$ solutions and the $L^1$ solutions of BSDEs. Then, we prove an existence result together with two comparison results for the maximal and minimal bounded solutions in Section 3, and two comparison theorems together with an existence and uniqueness result for the bounded solutions in Section 4. Finally, in Sections 5 and 6, we establish several existence results, comparison theorems, and existence and uniqueness results for the (maximal and minimal) $L^p\ (p> 1)$ and $L^1$ solutions of BSDEs respectively.

\section{Comparison theorems of $L^p\ (p>1)$ solutions and $L^1$ solutions}

In this section, we will establish several comparison theorems for the (maximal and minimal) $L^p\ (p> 1)$ solutions and the $L^1$ solutions of BSDEs. Let us first introduce the following assumptions on the generator $g$, where we assume that $0<T\leq +\infty$:\vspace{0.2cm}

{\bf (2A1)} There exists a function $u(\cdot)\in \US$ and a nondecreasing continuous function $\rho(\cdot):\R_+\mapsto \R_+$ with linear growth such that $\as$,
$$
(g(\omega,t,y_1,z)-g(\omega,t,y_2,z))\ {\rm sgn}(y_1-y_2) \leq u(t)\rho(|y_1-y_2|),\ \ \RE\ y_1,y_2,z.
$$
Assume further that $\rho(0)=0$, $\rho(u)>0$ for $u>0$ and $\int_{0^+}{1\over \rho(u)}\ {\rm d}u=+\infty$.\vspace{0.2cm}

{\bf (2A2)} There exists a function $v(\cdot)\in \VS$ and a nondecreasing continuous function $\phi(\cdot):\R_+\mapsto \R_+$ with $\phi (0)=0$ such that $\as$,
$$|g(\omega,t,y,z_1)-g(\omega,t,y,z_2)|\leq v(t)\phi(|z_1-z_2|),\ \ \RE\ y,z_1,z_2.$$
Without loss of generality, here and henceforth we can always assume that for all $x\in \R_+$, $0\leq\phi(x)\leq ax+b$. Furthermore, we assume also that $v(t)\in \US$ in the case where $b\neq 0$.\vspace{0.2cm}

{\bf (2A3)} $\as$, $g(\omega,t,\cdot,\cdot):\R\times \R^d\longmapsto \R$ is continuous.\vspace{0.2cm}

{\bf (2A4)} There exist two functions $u(\cdot)\in \US$, $v(\cdot)\in \VS$ and a process $f_t\in L^2(\Omega;L^1(\T;\R_+))$ such that $\as$,
$$
|g(\omega,t,y,z)| \leq f_t(\omega)+u(t)|y|+v(t)|z|, \ \ \RE\ y,z.
$$

{\bf (2A5)} There exists a constant $\alpha\in (0,1)$, a function $\lambda (\cdot):\T\mapsto \R_+$ and an $(\F_t)$-progressively measurable nonnegative processes $(f_t)_{t\in\T}\in L^1(\T\tim\Omega)$ such that $\as$,
$$
|g(\omega,t,y,z)-g(\omega,t,y,0)|\leq \lambda (t)(f_t(\omega)+|y|+|z|)^\alpha,\ \  \RE\ y, z.
$$
We also assume that
$$\int_0^T(\lambda (t)+\lambda^{1\over 1-\alpha}(t)+\lambda^{2\over 2-\alpha}(t))\ {\rm d}t<+\infty.$$

{\bf (2A5')} There exists a constant $\alpha\in (0,1)$ and a function $\lambda (\cdot):\T\mapsto \R_+$ such that $\as$,
$$
|g(\omega,t,y,z)-g(\omega,t,y,0)|\leq \lambda (t)|z|^\alpha,\ \  \RE\ y, z.
$$
We also assume that
$$\int_0^T \lambda^{2\over 2-\alpha}(t)\ {\rm d}t<+\infty.$$

{\bf Remark 2.1} In (2A1), we do not assume that $\rho(\cdot)$ is a concave function required by Theorem 2 of \citet{FJT11}.\vspace{0.1cm}

The following lemma will play an important role in the proof of main results of this paper.\vspace{0.1cm}

{\bf Lemma 2.1}\ Assume that $0<T\leq +\infty$, $\{b_n\}_{n=1}^{+\infty}$ is a nonnegative and non-increasing real sequence, $\beta(t)\in \US$ and $\psi(\cdot):\R_+\mapsto \R_+$ is a nondecreasing continuous function with linear growth. Let $\{(u_n(t))_{t\in \T}\}_{n=1}^{+\infty}$ be a sequence of non-negative $(\F_t)$-progressively processes satisfying $$\E_n\left[\sup\limits_{t\in [0,T]}u_n(t)\right]<+\infty$$
and for each $t\in\T$,
\begin{equation}
u_n(t)\leq b_n+\E_n\left[\left.\int_t^T \beta(s)\psi(u_n(s))\ {\rm d}s\right|\F_t\right]\ \ps,
\end{equation}
where $\E_n[X|\F_t]$ represents the conditional expectation of the random variable $X$ with respect to $\F_t$ under a probability measure $P_n$ which is defined on $(\Omega,\F_T)$ and may depend on $n$. If $\Lim b_n=0$, $\psi(0)=0$, $\psi(u)>0$ for $u>0$, and
\begin{equation}
\int_{0^+} {1\over \psi(u)}\ {\rm d}u=+\infty,
\end{equation}
then for each $t\in \T$,
\begin{equation}
\Lim u_n(t)=0\ \ \ps.\vspace{0.2cm}
\end{equation}

{\bf Proof.}\ Since $\psi$ is of linear growth, we can get the existence of a constant $k$ such that $\psi(x)\leq k(1+x)$ for all $x\in \R_+$. Then, in view of $b_n\leq b_1$, $\beta(\cdot)\in \US$ and the Fubinin Theorem, by (2) we can obtain that for each $t\in \T$, $\ps$,
$$\E_n[u_n(r)|\F_t]\leq b_1+k\int_0^T\beta (s)\ {\rm d}s+k\int_r^T \beta(s)\E_n[u_n(s)|\F_t]\ {\rm d}s,\ \ r\in [t,T].$$
Thus, Gronwall's inequality yields that for each $t\in \T$, $\ps$,
$$
\E_n[u_n(r)|\F_t]\leq  \left(b_1+k\int_0^T\beta (s)\ {\rm d}s\right)e^{k\int_r^T \beta (s)\ {\rm d}s},\ \ r\in [t,T].
$$
Letting $r=t$ in the above inequality we get that
for each $t\in \T$,
\begin{equation}
\sup\limits_{n\geq 1}u_n(t)\leq C:=\left(b_1+k\int_0^T\beta (s)\ {\rm d}s\right)e^{k\int_0^T \beta (s)\ {\rm d}s}\ \ \ \ps.
\end{equation}

In the sequel, in view of the linear growth of $\psi$, for each $n\geq 1$ we can define the function $\psi_n:\R_+\To \R_+$ as follows
$$\psi_n(x)=\sup_{y\in\R_+}\{\psi(y)-(n+2k)|x-y|\}.\vspace{-0.2cm}$$
It is well known that $\psi_n$ is well
defined and Lipschitz. Moreover, the sequence
$\{\psi_n\}_{n=1}^{+\infty}$ is non-increasing and converges to
$\psi$. Thus, for each $n\geq 1$, noticing that $\beta(\cdot)\in \US$, we can let $v_n:\R_+\To \R_+$ be the solution of the following backward ordinary differential equation (ODE for short):
$$v_n(t)=b_n+\int_t^T \beta (s) \psi_n
(v_n(s)){\rm d}s, \ \ t\in \T.
$$
Since
$\{\psi_n\}_{n=1}^\infty$ and $\{b_n\}_{n=1}^\infty$ are both non-increasing sequences, we know that $v_{n+1}\leq v_{n}$ for
each $n\geq 1$. This implies that, noticing that $\Lim b_n=0$ and that $\{\psi_n\}_{n=1}^{+\infty}$ converges to $\psi$ as $n\To \infty$, the sequence $\{v_n\}_{n=1}^\infty$ converges pointwisely to a function $v:\R_+\To \R_+$ which satisfies
$$v(t)=\int_t^T \beta (s) \psi (v(s)){\rm d}s, \ \ t\in \T.$$
In view of (3) and the fact that $\beta(\cdot)\in \US$, Bihari's inequality (see \citet{Bih56} for details) yields that $v(t)=0$ for each $t\in \T$.

Now for $n,j\geq 1$, let $v_n^j$ be the function defined
recursively as follows:
$$v^1_n(t)\equiv C,$$
\begin{equation}
v^{j+1}_n(t)=b_n+\int_t^T \beta (s) \psi_n (
v^j_n(s)){\rm d}s, \ j\geq 1,\ \ \ t\in \T,\vspace{0.1cm}
\end{equation}
where $C$ is defined in (5). Since $\psi_n$ is Lipschitz and $\beta(\cdot)\in \US$, we know that $v^j_n\To v_n$ as $j\To \infty$. On the other hand, it is easily seen by induction that for all $n,j\geq 1$ and each $t\in \T$,
\begin{equation}
u_n(t)\leq v^j_n(t)\ \ \ \ps.\vspace{0.1cm}
\end{equation}
Indeed, for $j=1$ the formula holds true by (5). Suppose it also holds for some $j$, then
$$\psi(u_n(s))\leq \psi(v_n^j(s))\leq
\psi_n(v_n^j(s)),\ \ \ s\in \T.$$
In view of (2), the previous inequality and (6), we can deduce that for all $n\geq 1$ and each $t\in \T$,
$$u_n(t)\leq v_n^{j+1}(t)\ \ \ \ps.$$
Thus, (7) follows. Finally, taking the limit in (7) as first
$j\To \infty$, and then $n\To \infty$, we
obtain (4). The proof of Lemma 2.1 is then completed.\vspace{0.2cm}\hfill $\Box$

The following Theorem 2.1 establishes a general comparison theorem for $L^p\ (p>1)$ solutions of BSDEs. It can be regarded as a generalization of Theorem 2 in \citet{FJT11}, where only are the $L^2$ solutions considered, and the concavity condition of $\rho(\cdot)$ in (2A1) is also required.\vspace{0.1cm}

{\bf Theorem 2.1} \ Assume that $0<T\leq +\infty$, $g$ and $g'$ are two generators of BSDEs and  $(y_\cdot,z_\cdot)$ and $(y'_\cdot,z'_\cdot)$ are, respectively, a solution of BSDE$(\xi,g)$ and BSDE$(\xi',g')$. Assume further that $(y_\cdot-y'_\cdot)^+\in \s$. If $\xi\leq \xi'\ \ps$ and one of the following two statements holds true:

(i)\ \ $g$ satisfies (2A1) and (2A2), and
\begin{equation}
\mathbbm{1}_{ y_t>y'_t}(g(t,y'_t,z'_t)-g'(t,y'_t,z'_t))\leq 0\ \ \as;
\end{equation}

(ii) $g'$ satisfies (2A1) and (2A2), and
\begin{equation}
\mathbbm{1}_{ y_t> y'_t}(g(t,y_t,z_t)-g'(t,y_t,z_t))\leq 0 \ \ \as,\vspace{0.1cm}
\end{equation}

\noindent then for each $t\in\T$, we have\vspace{-0.1cm}
$$y_t\leq y'_t\ \ \ps.\vspace{-0.1cm}$$

{\bf Proof.}\ \ We will only prove the case (i). Another case can be proved in the same way. Let us fix $k\in {\N}^*$ and
denote $\tau_k$ the stopping time
$$\tau_k=\inf \left\{t\in\T: \int_0^t\left(|z_s|^2+|z'_s|^2
\right)\ {\rm d}s\geq k\right\}\wedge T.$$
Tanaka's formula leads to the equation, setting
$\hat{y}_t=y_t-y'_t,\ \hat{z}_t=z_t-z'_t$,
\begin{equation}
\Dis \hat{y}_{t\wedge \tau_k}^+\leq \hat{y}_{\tau_k}^+
+\int_{t\wedge \tau_k}^{\tau_k}
 \mathbbm{1}_{ \hat{y}_s> 0}(g(s,y_s,z_s)-g'(s,y'_s,z'_s))
 \ {\rm d}s-\int_{t\wedge \tau_k}^{\tau_k}\mathbbm{1}_{ \hat{y}_s> 0}
\hat{z}_s\cdot {\rm d}B_s.
\end{equation}
First of all, since $\mathbbm{1}_{ \hat{y}_s> 0}(g(s,y'_s,z'_s)-g'(s,y'_s,z'_s))$ is
non-positive, we have
$$
\begin{array}{lll}
&&\mathbbm{1}_{ \hat{y}_s> 0}(g(s,y_s,z_s)-g'(s,y'_s,z'_s))\\
&=&\mathbbm{1}_{ \hat{y}_s> 0}(g(s,y_s,z_s)-g(s,y'_s,z'_s))
+\mathbbm{1}_{ \hat{y}_s> 0}(g(s,y'_s,z'_s)-g'(s,y'_s,z'_s))\\
&\leq & \mathbbm{1}_{ \hat{y}_s> 0}(g(s,y_s,z_s)-g(s,y'_s,z_s))+\mathbbm{1}_{ \hat{y}_s> 0}(g(s,y'_s,z_s)
-g(s,y'_s,z'_s))
\end{array}$$
and we deduce, using assumptions (2A1) and (2A2) for $g$, that
\begin{equation}
\mathbbm{1}_{ \hat{y}_s> 0}(g(s,y_s,z_s)-g'(s,y'_s,z'_s))\leq
u(s)\rho(\hat{y}_s^+)+\mathbbm{1}_{ \hat{y}_s> 0}v(s)\phi(|\hat{z}_s|).
\end{equation}
From the proof of Theorem 1 in \citet{FJD10} we know that, with $c=a+b$,
\begin{equation}
\phi(x)\leq (n+2c)x+{\mathbf 1}_{ b\neq 0}\phi\left({2c\over n+2c}\right),\ \ \RE\ x\in \R_+,\ \RE\ n\geq 1.
\end{equation}
Combining (10)-(12) yields that for each $n\geq 1$ and each $t\in\T$,
\begin{equation}
\begin{array}{lll}
\Dis \hat{y}_{t\wedge \tau_k}^+ &\leq &\Dis  a_n +\hat{y}_{\tau_k}^+ +\int_{t\wedge \tau_k}^{\tau_k}
\left[u(s)\rho(\hat{y}_s^+)+\mathbbm{1}_{ \hat{y}_s>
0}(n+2c)v(s)|\hat{z}_s|\right]{\rm d}s\\
&&\Dis -\int_{t\wedge \tau_k}^{\tau_k}\mathbbm{1}_{ \hat{y}_s> 0}\hat{z}_s\cdot {\rm d}B_s\\
&=&\Dis a_n+\hat{y}_{\tau_k}^+ +\int_{t\wedge \tau_k}^{\tau_k}u(s)\rho(\hat{y}_s^+)\ {\rm
d}s\\
&&\Dis -\int_{t\wedge \tau_k}^{\tau_k}\mathbbm{1}_{ \hat{y}_s> 0}\hat{z}_s\cdot [-{(n+2c)v(s)\hat{z}_s\over
|\hat{z}_s|}\mathbbm{1}_{ |\hat{z}_s|\neq 0}\ {\rm d}s+{\rm d}B_s],
\end{array}
\end{equation}
where
\begin{equation}
a_n={\mathbf 1}_{ b\neq 0}\phi({2c\over n+2c})\cdot\int_0^T v(s)\ {\rm d}s\To 0\ \ {\rm by\ (2A2)\ as}\ n\To \infty.\vspace{0.2cm}
\end{equation}
Let
$P_n$ be the probability on $(\Omega,\F_T)$ which is equivalent to $P$ and defined by
$${{\rm d}P_n\over {\rm d}P}:=\exp \left\{(n+2c)\int_0^T {v(s)\hat{z}_s\over
|\hat{z}_s|}\mathbbm{1}_{ |\hat{z}_s|\neq 0}\cdot {\rm d}B_s-{1\over
2}(n+2c)^2\int_0^T \mathbbm{1}_{ |\hat{z}_s|\neq 0}v^2(s)\ {\rm d}s\right\}.$$
It is worth noticing that ${\rm d}P_n/{\rm d}P$ has moments of all order since $v(\cdot)\in \VS$. By Girsanov's theorem, under $P_n$ the process
$$B_n(t)=B_t-\int_0^t {(n+2c)v(s)\hat{z}_s\over
|\hat{z}_s|}\mathbbm{1}_{ |\hat{z}_s|\neq 0}\ {\rm d}s,\ \
t\in\T,\vspace{-0.1cm}$$
is an $(\F_t,P_n)-$Brownian motion.
Moreover, the process
$$\left(\int_0^{t\wedge\tau_k} \mathbbm{1}_{ \hat{y}_s> 0}\hat{z}_s\cdot\ {\rm d}B_n(s)\right)_{0\leq t\leq T}$$
is an $(\F_t,P_n)-$martingale. Let ${\E}_n[X|\F_t]$ represent the conditional expectation of the random variable $X$ with respect to $\F_t$ under $P_n$. By taking the conditional
expectation with respect to $\F_t$ under $P_n$ in (13) we get that for each $n\geq 1$ and $t\in\T$,
$$
\Dis \hat{y}_{t\wedge\tau_k}^+\leq a_n+
\E_n\left[\left.\hat{y}_{\tau_k}^+\right|\F_{t}\right]
+\E_n\left[\left.\int_{t\wedge\tau_k}^{\tau_k}
u(s)\rho\left(\hat{y}_s^+\right){\rm d}s\right|\F_{t}\right].
$$
Furthermore, in view of the facts that $\tau_k\To T$ as $k\To \infty$, $(\hat y_\cdot)^+$ belongs to $\s$, $\xi\leq \xi'$ and $u(\cdot)\in \US$, letting $k\To \infty$ in the above inequality and using Lebesgue's dominated convergence theorem yields that for each $t\in\T$,
$$
\Dis \hat{y}_{t}^+\leq a_n+\E_n\left[\left.\int_{t}^{T}u(s)\rho\left(\hat{y}_s^+ \right){\rm d}s\right|\F_{t}\right].
$$
Thus, in view of (14), applying Lemma 2.1 with
$u_n(t)\equiv \hat{y}_t^+$,
$b_n=a_n$, $\beta (s)=u(s)$ and $\psi(u)=\rho(u)$ yields that for each $t\in \T$,
$$\hat{y}_t^+=\Lim \hat{y}_t^+=\Lim u_n(t)=0\ \ \ \ps.$$
That is to say, for each $t\in\T$, $y_t\leq
y'_t\ \ps$. Theorem 2.1 is then proved. \vspace{0.2cm}\hfill $\Box$

It should be especially noted that the presence of the indicate function makes that (8) and (9) can be more easily satisfied than the usual form. This important observation is one of the start points of this paper.\vspace{0.2cm}

The following Remark 2.2 gives a easily verifiable condition to ensure that (8) or (9) holds true, which will be used several times and play an important role  later.\vspace{0.1cm}

{\bf Remark 2.2}\ Assume that $c$ is a constant. If
$$\ps,\ \RE\ t\in\T,\ y_t\leq c$$
and
$$\as,\ \RE\ y<c,\ \RE\ z\in\R^d,\ g(t,y,z)\leq g'(t,y,z),\vspace{0.1cm}$$
then (8) holds true. Similarly, if
$$\ps,\ \RE\ t\in\T,\ y'_t\geq c$$
and
$$\as,\ \RE\ y>c,\ \RE\ z\in\R^d,\ g(t,y,z)\leq g'(t,y,z),$$
then (9) holds true.\vspace{0.2cm}

Let us further introduce the following Lemma 2.2, which comes from Theorem 1 in \citet{FJT11}. It will be generalized in Section 5.\vspace{0.1cm}

{\bf Lemma 2.2}\ Assume that $0<T\leq +\infty$ and the generator $g$ satisfies (2A3) and (2A4). Then for each $\xi\in L^2(\Omega,\F_T,P)$, BSDE$(\xi,g)$ has both a minimal $L^2$ solution and a maximal $L^2$ solution.\vspace{0.2cm}

The following Theorem 2.2 establishes a general comparison theorem for the maximal $L^2$ solutions, which will be also generalized in Section 5.\vspace{0.1cm}

{\bf Theorem 2.2}\ \ Assume that $0<T\leq +\infty$, $g$ and $g'$ are two generators of BSDEs, and $(y_\cdot,z_\cdot)$ is any $L^2$ solution of BSDE$(\xi,g)$. Assume further that $g'$ satisfies (2A3) and (2A4), and $(y'_\cdot,z'_\cdot)$ is the maximal $L^2$ solution of BSDE$(\xi',g')$ by Lemma 2.2. If $\xi\leq \xi'\ \ps$ and (9) holds true, then for each $t\in\T$,
$$y_t\leq y'_t\ \ \ps.$$

{\bf Proof.}\ \ Since $g'$ satisfies (2A3) and (2A4), by the proof of Theorem 1 in \citet{FJT11} we know that if we take for each $n\geq 1$ and $(\omega,t,y,z)\in \Omega\times\T\times\R\times \R^d$,
$$
{}^ng'(\omega,t,y,z):=\sup\limits_{(u,v)\in \R^{1+d}}
\{g'(\omega,t,u,v)-nu(t)|y-u|-nv(t)|z-v|\},
$$
then ${}^n g'(t,\cdot,\cdot)\downarrow g'(t,\cdot,\cdot)$ and ${}^ng'(t,\cdot,\cdot)$ satisfies
\begin{equation}
|{}^ng'(t,y_1,z_1)-{}^ng'(t,y_2,z_2)|\leq
nu(t)|y_1-y_2|+nv(t)|z_1-z_2|,\ \ \RE\ y_1,y_2,z_1,z_2.
\end{equation}
Furthermore, from the proof of Theorem 1 in \citet{FJT11}, we also know that $({}^ny'_\cdot,{}^nz'_\cdot)$, the unique $L^2$ solution of BSDE$(\xi',{}^ng')$, satisfies
\begin{equation}
{}^ny'_\cdot\downarrow y'_\cdot.
\end{equation}
Combining (9) and (16) as well as the fact ${}^n g'\geq g'$ we can get that for each $n\geq 1$,
$$
\mathbbm{1}_{ y_t>{}^ny'_t}g(t,y_t,z_t)\leq \mathbbm{1}_{ y_t>{}^ny'_t}g'(t,y_t,z_t)\leq
\mathbbm{1}_{ y_t>{}^ny'_t}{}^ng'(t,y_t,z_t)\ \ \as.
$$
Thus, in view of (15) and the above inequality, Theorem 2.1 yields that
\begin{equation}
\RE\ n\geq 1,\ \RE\ t\in \T,\ \ y_t\leq {}^ny'_t\ \ \ps.
\end{equation}
In view of (16), the conclusion follows by letting $n\To \infty$ in (17). \vspace{0.2cm}\hfill $\Box$

By similar argument to Theorem 2.2 we can get the following Theorem 2.3.\vspace{0.1cm}

{\bf Theorem 2.3}\ Assume that $0<T\leq +\infty$, $g$ and $g'$ are two generators of BSDEs, and $(y'_\cdot,z'_\cdot)$ is any $L^2$ solution of BSDE$(\xi',g')$. Assume further that $g$ satisfies (2A3) and (2A4), and $(y_\cdot,z_\cdot)$ is the minimal $L^2$ solution of BSDE$(\xi,g)$ by Lemma 2.2. If $\xi\leq \xi',\ \ps$ and (8) holds true, then for each $t\in\T$,
$$y_t\leq y'_t\ \ \ \ps.$$

{\bf Corollary 2.1} Assume that $0<T\leq +\infty$ and both $g$ and $g'$ satisfy (2A3) and (2A4). Let $(y_\cdot,z_\cdot)$ and $(y'_\cdot,z'_\cdot)$ be, respectively, the maximal (resp. minimal) $L^2$ solution of BSDE$(\xi,g)$ and BSDE$(\xi',g')$ by Lemma 2.2. If $\xi\leq \xi'\ \ps$ and for each $(y,z)\in \R\times \R^d$,
$$g(t,y,z)\leq g'(t,y,z)\ \ \ \as,$$
then for each $t\in \T$,
$$y_t\leq y'_t\ \ \ \ps.$$

The following Theorem 2.4 establishes a general comparison theorem for $L^1$ solutions of BSDEs, which improves virtually several corresponding comparison results obtained respectively in \citet{Bri06}, \citet{FL10}, \citet{Xiao12}, \citet{Fan12} and \citet{Tian13} even for the case of the finite time horizon.\vspace{0.1cm}

{\bf Theorem 2.4}\ \ Assume that $0<T\leq +\infty$, $g$ and $g'$ are two generators of BSDEs, and $(y_\cdot,z_\cdot)$ and $(y'_\cdot,z'_\cdot)$ are, respectively, an $L^1$ solution of BSDE$(\xi,g)$ and BSDE$(\xi',g')$. If $\xi\leq \xi'\ \ \ps$ and one of the following two statements is satisfied:\vspace{0.2cm}

(i)\ \ $g$ satisfies (2A1), (2A2) and (2A5) (or (2A5')), and (8) holds true;\vspace{0.1cm}

(ii) $g'$ satisfies (2A1), (2A2) and (2A5) (or (2A5')), and (9) holds true,\vspace{0.2cm}

\noindent then for each $t\in\T$,
$$y_t\leq y'_t\ \ \ \ps.$$

{\bf Proof.}\ It follows from Theorem 2.1 that we need only to show that $(y_\cdot-y'_\cdot)^+$ belongs to $\s$ under the assumptions of Theorem 2.4.

We now assume that $\xi\leq \xi'\ \ \ps$, $g$ satisfies (2A1), (2A2) and (2A5), and (8) holds true. The same arguments as follows can prove the other cases. Let us fix $k\in \N^*$ and
denote the stopping time
$$\tau_k:=\inf \left\{t\in\T: \int_0^t\left(|z_s|^2+|z'_s|^2
\right)\ {\rm d}s\geq k\right\}\wedge T.
$$
Tanaka's formula leads to the equation, setting $\hat{y}_t=y_t-y'_t,\
\hat{z}_t=z_t-z'_t$,
$$
\Dis \hat{y}_{t\wedge \tau_k}^+\leq \hat{y}_{\tau_k}^+
+\int_{t\wedge \tau_k}^{\tau_k}
 \mathbbm{1}_{ \hat{y}_s> 0}(g(s,y_s,z_s)-g'(s,y'_s,z'_s))
 \ {\rm d}s-\int_{t\wedge \tau_k}^{\tau_k}\mathbbm{1}_{ \hat{y}_s> 0}
\hat{z}_s\cdot {\rm d}B_s.$$
Since $\mathbbm{1}_{ \hat{y}_s> 0}(g(s,y'_s,z'_s)-g'(s,y'_s,z'_s))$ is
non-positive, we have
$$
\begin{array}{lll}
&&\mathbbm{1}_{ \hat{y}_s> 0}(g(s,y_s,z_s)-g'(s,y'_s,z'_s))\\
&=&\mathbbm{1}_{ \hat{y}_s> 0}(g(s,y_s,z_s)-g(s,y'_s,z'_s))
+\mathbbm{1}_{ \hat{y}_s> 0}(g(s,y'_s,z'_s)-g'(s,y'_s,z'_s))\\
&\leq & \mathbbm{1}_{ \hat{y}_s> 0}(g(s,y_s,z_s)-g(s,y'_s,z_s))+\mathbbm{1}_{ \hat{y}_s> 0}(g(s,y'_s,z_s)
-g(s,y'_s,z'_s))
\end{array}$$
and we deduce, using assumptions (2A1) and (2A5) of $g$, that
$$
\mathbbm{1}_{ \hat{y}_s> 0}(g(s,y_s,z_s)-g'(s,y'_s,z'_s))\leq
u(s)\rho(\hat{y}_s^+)+2\lambda(s)(f_s+|y'_s|+|z_s|+|z'_s|)^{\alpha}.
$$
Thus, we get that, with $\phi_s:=2\lambda(s)(f_s+|y'_s|+|z_s|+|z'_s|)^{\alpha}$,
$$
\Dis \hat{y}_{t\wedge \tau_k}^ +\leq \hat{y}_{\tau_k}^+
+\int_{t\wedge \tau_k}^ { \tau_k}\left(u(s)\rho(\hat{y}_s^+)+\phi_s\right){\rm d}s-\int_{t\wedge
\tau_k}^ { \tau_k}\mathbbm{1}_{ \hat{y}_s> 0}\hat{z}_s\cdot {\rm d}B_s,
$$
and then that
\begin{equation}
\Dis\hat{y}_{t\wedge \tau_k}^ + \leq
\Dis\E\left[\hat{y}_{\tau_k}^+ +\left.\int_{t\wedge
\tau_k}^{\tau_k}
\left(u(s)\rho(\hat{y}_s^+)+\phi_s\right) {\rm d}s\right|\F_{t}\right].\vspace{0.1cm}
\end{equation}

Furthermore, since $\rho(\cdot)$ is of linear growth, we can find a pair of positive constants $k_1$ and $k_2$ such that
\begin{equation}
\rho(u)\leq k_1+k_2u,\ \ \RE\  u\geq 0.\vspace{0.1cm}
\end{equation}
Then, since both $(y_\cdot,z_\cdot)$ and $(y'_\cdot,z'_\cdot)$ are the $L^1$ solutions, we can send $k$ to $\infty$ in (18) and use Lebesgue's dominated convergence theorem and Fubini's theorem, in view of $\xi\leq \xi'$,  $\tau_k\To T$ as $k\To \infty$, $u(\cdot)\in \US$ and (19), to get that, for each $t\in \T$,
$$
\begin{array}{lll}
\Dis \hat{y}_t^+ &\leq & \Dis \E\left[\left.\int_0^T
\phi_s\ {\rm
d}s\right|\F_t\right]+\E\left[\left.\int_t^T u(s)\rho(\hat{y}_s^+){\rm
d}s\right|\F_t\right]\vspace{0.2cm}\\
&\leq &\Dis k_1\int_0^T u(s)\ {\rm d}s+\E\left[\left.\int_0^T
\phi_s\ {\rm d}s\right|\F_t\right]+k_2\int_t^T u(s)
\E\left[\left.\hat{y}_s^+\right|\F_t\right]{\rm d}s,
\end{array}\vspace{0.2cm}$$
and then for each $r\in [t,T]$,
$$\E\left[\left.\hat{y}_r^+\right|\F_t\right]\leq k_1\int_0^T u(s)\ {\rm d}s+
\E\left[\left.\int_0^T \phi_s\ {\rm d}s\right|\F_t\right]+k_2\int_r^T u(s)\E\left[\left.\hat{y}_s^+\right|\F_t\right]{\rm d}s.\vspace{0.2cm}$$
Gronwall's inequality yields that for each $r\in [t,T]$,
$$\E\left[\left.\hat{y}_r^+\right|\F_t\right]\leq
\left(k_1\int_0^T u(s)\ {\rm d}s+\E\left[\left.\int_0^T \phi_s\ {\rm
d}s\right|\F_t\right]\right)\cdot e^{k_2\int_r^T u(s){\rm d}s},$$ from which, by letting $r=t$, we have
\begin{equation}
\hat{y}_t^+\leq \left(k_1\int_0^T u(s)\ {\rm d}s+\E\left[\left.\int_0^T
\phi_s\ {\rm d}s\right|\F_t\right]\right)\cdot e^{k_2\int_0^T u(s) {\rm d}s}.
\end{equation}

Now, let $\beta$ be any constant which belongs to $(\alpha,1)$. Then we have
\begin{equation}
\E\left[\left(\int_0^T
\phi_s\ {\rm d}s\right)^{\beta\over \alpha}\right]<+\infty.
\end{equation}
Indeed, H\"{o}lder's inequality yields that
$$\int_0^T \lambda(s)f^\alpha_s\ {\rm d}s\leq\left(\int_0^T \lambda^{1\over 1-\alpha}(s)\ {\rm d}s\right)^{1-\alpha}\left(
\int_0^T f_s\ {\rm d}s\right)^\alpha,
$$
$$
\int_0^T \lambda(s)|z_s|^\alpha\ {\rm d}s\leq \left(\int_0^T\lambda^{2\over 2-\alpha}(s)\ {\rm d}s\right)^{2-\alpha\over 2}\left(\int_0^T|z_s|^2\ {\rm d} s\right)^{\alpha\over 2},\vspace{0.1cm}
$$
and $z'_s$ has a similar estimate. Besides,
$$
\int_0^T\lambda(s)|y'_s|^\alpha\ {\rm d}s\leq\int_0^T\lambda(s)\ {\rm d}s\cdot\sup_{t\in\T}|y'_t|^{\alpha}.
$$
Thus, coming back to the definition of $\phi_s$ and noticing the assumptions of the deterministic function $\lambda(\cdot)$ and the facts that $f_s\in L^1(\Omega\times \T)$, $(z_t)_{t\in\T}$ and $(z'_t)_{t\in\T}$ belong to the space $\M^\beta$, and $(y'_t)_{t\in\T}$ belongs to the space $\s^\beta$, we can obtain (21).\vspace{0.1cm}

Finally, taking the supremum over $t$ and the mathematical expectation after taking the power of $\beta \over \alpha$ in both sides of (20) and then making use of Doob's inequality, we can get that there exists a constant $\bar k>0$ such that, in view of (21),
$$
\E[\sup_{t\in\T}|\hat{y}_t^+|^{\beta\over \alpha}]\leq \bar k+\bar k\E\left[\left(\int_0^T
\phi_s\ {\rm d}s\right)^{\beta\over \alpha}\right]<+\infty,
$$
which means that $(y_\cdot-y'_\cdot)^+\in \s$. Thus, we complete the proof of Theorem 2.4. \hfill $\Box$\vspace{0.2cm}

\section{Existence of bounded solutions}

In this section, we will put forward and prove an existence result and two comparison results for the maximal and minimal bounded solutions. We denote by $\LS$ the set of continuous and strictly positive functions $l(x):\R\mapsto \R^+$ satisfying
$$
\int_{-\infty}^{0}{{\rm d}x\over l(x)}=\int_0^{+\infty}{{\rm d}x\over l(x)}=+\infty.
$$
We will use the following assumptions on the generator $g$, where $0<T\leq +\infty$:\vspace{0.2cm}

{\bf (3A1)} There exists a function $u(\cdot)\in \US$ and two nonnegative continuous functions $l(\cdot)\in\LS$ and $h(\cdot):\R\mapsto \R_+$ such that $\as$,
$$g(\omega,t,y,z)\ {\rm sgn}(y)\leq u(t)l(y)+h(y)|z|^2,\ \ \RE\ y,z.$$

{\bf (3A2)} There exists a function $\bar u(\cdot)\in \US$ and two nonnegative continuous functions $\bar \varphi(\cdot),\ \bar h(\cdot):\R\mapsto \R_+$ such that $\as$,
$$|g(\omega,t,y,z)|\leq \bar u(t)\bar \varphi(y)+\bar h(y)|z|^2,\ \ \RE\ y,z.$$

To state a main result of this section, we introduce the following lemma 3.1.\vspace{0.1cm}

{\bf Lemma 3.1}\ Assume that $0<T\leq +\infty$, $u(\cdot)\in \US$ and $l(x):\R\mapsto \R^+$ is a continuous function. Then $l(\cdot)\in \LS$ if and only if for each $-\infty<a\leq 0\leq b<+\infty$, the following two backward ODEs:
\begin{equation}
L_t=a-\int_t^T u(s)\ l(L_s)\ {\rm d}s
\end{equation}
and
\begin{equation}
U_t=b+\int_t^T u(s)\ l(U_s)\ {\rm d}s\vspace{0.2cm}
\end{equation}
have both global bounded solutions on $\T$.

Moreover, if $l(\cdot)\in \LS$, then (22) and (23) have unique global bounded solutions $L_t$ and $U_t$ on $\T$, and for each $t\in \T$, we have
$$L_0\leq L_t\leq a\leq 0\leq b\leq U_t\leq U_0.$$

{\bf Proof.}\ In the proof of Lemma 1 in \citet{Lep98}, by replacing the term $T-t$ with $$\int_t^T u(s)\ {\rm d}s$$
we can complete the proof of this lemma.\vspace{0.2cm}\hfill $\Box$

The following Theorem 3.1 is one of main results in this section, which generalizes virtually the corresponding existence results for bounded solutions of BSDEs obtained respectively in \citet{Kob00}, \citet{Lep98}, \citet{Bri08} and \citet{Bri07} even for the case of the finite time horizon.\vspace{0.1cm}

{\bf Theorem 3.1} Assume that $0<T\leq +\infty$ and that $g$ satisfies (2A3), (3A1) and (3A2). Then for each $\xi\in L^\infty(\Omega,\F_T,P)$, BSDE$(\xi,g)$ has both a minimal one and a maximal one among all bounded solutions $(Y,Z)$. Moreover, for each $t\in \T$,
$$L_0\leq L_t\leq Y_t\leq U_t\leq U_0\ \ \ps,$$
where $(L,U)$ are the unique solutions of (22) and (23) with $a=-\|\xi\|_{\infty}$ and $b=\|\xi\|_{\infty}$.\vspace{0.2cm}

To prove Theorem 3.1, we need the following three lemmas. First, by similar argument to the proof of Lemma 3 in \citet{Lep98} and in view of Lemma 2.2 and Theorem 2.1 (or Theorem 1.2 in \citet{Chen00}), we can obtain the following Lemma 3.2.\vspace{0.1cm}

{\bf Lemma 3.2}\ Assume that $0<T\leq +\infty$, $u(\cdot)\in \US$, $v(\cdot)\in \VS$, $f_1(\cdot):\R\mapsto \R$ is a continuous function with linear growth and the random function $f_2(\omega,t,z):\Omega\times\T\times\R^d\mapsto \R$ is $\mathcal P\times \mathcal B (\R^d)$-measurable, with $\as$,
$$
f_2(\omega,t,0)=0\ \ {\rm and}\ \ |f_2(\omega,t,z_1)-f_2(\omega,t,z_2)|\leq v(t)|z_1-z_2|,\ \ \RE\ z_1,z_2\in\R^d.
$$
If the backward ODE
$$J_t=a+\int_t^T u(s)f_1(J_s)\ {\rm d}s,\ \ a\in \R$$
has a unique solution $J$ on $\T$, then the BSDE
$$Y_t=a+\int_t^T (u(s)f_1(Y_s)+f_2(s,Z_s))\ {\rm d}s-\int_t^T Z_s\cdot {\rm d}B_s,\ \ t\in \T$$
has a unique solution given by $Z\equiv 0$ and $Y\equiv J$.\vspace{0.1cm}

The following Lemma 3.3 is the first step to prove Theorem 3.1.\vspace{0.1cm}

{\bf Lemma 3.3} Assume that $0<T\leq +\infty$ and $\eta\in L^\infty(\Omega,\F_T,P)$ which satisfies $0\leq \alpha\leq \eta\leq \beta\ \ \ps$. We assume without loss of generality that $\alpha\leq 1$ and $\beta\geq 1$. Suppose that $G(\omega,t,y,z):\Omega\times\T\times\R\times\R^d\mapsto \R$ vanishes when $y\leq 0$ and verifies for a constant $k>0$ and a function $u(\cdot)\in \US$ the following restriction:
$$\as,\ \RE\ y> 0,\  \RE\ z\in \R^d,\ \ -u(t)y-k|z|^2\leq G(\omega,t,y,z)\leq u(t)y.$$
Also we assume that $G$ is $\mathcal{P}\otimes\mathcal{B}(\R^{d+1})$ measurable and $\as$, $G(\omega,t,\cdot,\cdot)$ is continuous. Then the BSDE
$$
Y_t=\eta+\int_t^T G(s,Y_s,Z_s)\ {\rm d}s-\int_t^T Z_s\cdot {\rm d}B_s,\ \ t\in\T
$$
has a maximal bounded solution $(\theta,\Gamma)$. Moreover, for each $t\in \T$, we have $ Q_0\leq Q_t\leq \theta_t\leq S_t\leq S_0\ \ \ps$, where
$$Q_t:=\alpha \exp\left(-\int_t^T u(s)\ {\rm d}s\right),\ S_t:=\beta \exp\left(\int_t^T u(s)\ {\rm d}s\right).$$

{\bf Proof.}\ We will follow those steps used in the proof of Theorem 2 in \citet{Lep98}. We take $\kappa_n:\R^d\mapsto \R$ a sequence of smooth functions such that
$$1)\ 0\leq \kappa_n\leq 1; \ \ \ 2)\ \kappa_n(z)=1\ {\rm if}\ |z|\leq n;\ \ \ 3)\ \kappa_n(z)=0\ {\rm if}\ |z|\geq n+1.$$
For each $(\omega,t,y,z)\in \Omega\times\T\times\R\times\R^d$ and each $n\geq 1$, define
$$ G_n(\omega,t,y,z):=u(t)y \mathbbm {1}_{y>0}(1-\kappa_n(z\cdot e^{{t\over 2}}))
+\kappa_n(z\cdot e^{{t\over 2}}) G(\omega,t,y,z).$$
It is immediately seen that $G_n\downarrow G$ and that $ G_n$ is a continuous function of $(y,z)$ which satisfies that $\as$,
\begin{equation}
| G_n (\omega,t,y,z)|\leq k(n+1)^2e^{-t}+u(t)|y|,\ \ \RE\ y,z.
\end{equation}
Then by Lemma 2.2 we have the existence of a maximal $L^2$ solution $(\theta^n,\Gamma^n)$ for the equation
$$\theta_t^n=\eta+\int_t^T  G_n(s,\theta_s^n,\Gamma_s^n)\ {\rm d}s-\int_t^T \Gamma_s^n\cdot {\rm d}B_s.$$
Since $\as$,
$$\lambda_n(t,y,z):=-u(t)|y|-k|z|^2\kappa_n(z\cdot e^{{t\over 2}})\leq  G_n(\omega,t,y,z)\leq u(t)|y|,\ \ \RE\ y,z,$$
noticing that both $u(t)|y|$ and $\lambda_n$ satisfy (2A1)-(2A4), from Theorem 2.1 we can get that for each $t\in \T$,
$$Q^n_t\leq \theta^n_t\leq \bar S_t\ \ \ps,$$
where $(Q^n,Z^n)$ and $(\bar S,\bar Z)$ are, respectively, the unique $L^2$ solution of BSDE$(\alpha,\lambda_n)$ and BSDE$(\beta,u(t)|y|)$ by Lemma 2.2 and Theorem 2.1. Furthermore, 
it follows from Lemma 3.2 that $(\bar S,\bar Z)=(S,0)$ and for each $n\geq 1$, $(Q^n,Z^n)=(Q,0)$.

Consequently, $\theta^n$ is a decreasing and bounded sequence (in view of (24) and Corollary 2.1), then we have the existence of $\theta$ such that for each $t\in \T,\ \theta^n_t\downarrow \theta_t\ \ps$ and
$$\E\left[\int_0^Tu(t)|\theta^n_t-\theta_t|^2\ {\rm d}t\right]\To 0,\ \ {\rm as}\ \ n\To \infty.\vspace{0.2cm}$$
Moreover $\theta$ satisfies that for each $t\in \T,\ Q_0\leq Q_t\leq \theta_t\leq S_t\leq S_0\ \ps$.

In the sequel, following closely the proof procedure of Theorem 2 in \citet{Lep98} and noticing that $\as$, for each $(y,z)\in [Q_0,S_0]\times \R^d$,
$$
\sup\limits_{n\geq 1}|G_n (\omega,t,y,z)|\leq u(t)S_0+k|z|^2,\vspace{-0.2cm}
$$
we can prove that $\Gamma^n$ has a convergent subsequence in ${\rm M}^2$. Thus, take $\Gamma$ any accumulation point in ${\rm M}^2$ of $\Gamma^n$ then it is easy to verify that $(\theta,\Gamma)$ is a solution of BSDE$(\eta,G)$ (for more details see the proof of Theorem 1 in \citet{FJT11}).

Finally, for any bounded solution $(\hat Y,\hat Z)$ of BSDE$(\eta,G)$, noticing that $G_n\downarrow G$, (24) and the fact that $(\theta^n,\Gamma^n)$ is the maximal $L^2$ solution of BSDE$(\eta, G^n)$, by Theorem 2.2 we can conclude that for each $t\in \T$ and each $n\geq 1$,
$$\hat Y_t\leq \theta^n_t\ \ \ps,$$
and then $\hat Y\leq \theta$. Thus, Lemma 3.3 is proved.\vspace{0.2cm}\hfill $\Box$

By virtue of Lemma 3.3 we can prove the following Lemma 3.4.\vspace{0.1cm}

{\bf Lemma 3.4} Assume that $0<T\leq +\infty$ and the generator $g$ satisfies (2A3). Assume further that there exists a constant $\gamma>0$ and a function $u(\cdot)\in \US$ such that $\as$,
\begin{equation}
|g(\omega,t,y,z)|\leq u(t)+{\gamma\over 2} |z|^2,\ \ \RE\ y,z.
\end{equation}
Then for each $\xi\in L^\infty(\Omega,\F_T,P)$, BSDE$(\xi,g)$ has both a maximal one and a minimal one among all bounded solutions $(y_\cdot,z_\cdot)$. Moreover, for each $t\in \T$, we have
$$|y_t|\leq \|\xi\|_{\infty}+\int_t^T u(s)\ {\rm d}s\ \ \ \ps.\vspace{-0.2cm}$$

{\bf Proof.}\ We will first prove the existence of the maximal bounded solution by following those steps used in the proof of Theorem 1 in \citet{Lep98}. Let
$$\eta:=e^{\gamma \xi}\in L^\infty(\Omega,\F_T,P),\ \alpha:=e^{-\gamma \|\xi\|_{\infty}},\ \ \beta:=e^{\gamma \|\xi\|_{\infty}},\vspace{0.2cm}$$
$$Q_t=\alpha\exp{\left(-\gamma\int_t^T u(s)\ {\rm d}s\right)},\ S_t=\beta\exp{\left(\gamma\int_t^T u(s)\ {\rm d}s\right)},\ \ t\in\T,\vspace{0.2cm}$$
and for each $(\omega,t,y,z)\in \Omega\times\T\times\R\times\R^d$, define
$$G(\omega,t,y,z):=\mathbbm{1}_{ y>0}\left(\gamma yg\left(\omega,t,{\ln y\over \gamma},{z\over \gamma y}\right)-{1\over 2}{|z|^2\over y}\right).$$
It then follows from (25) that $0<\alpha\leq 1\leq \beta$, $\alpha\leq \eta\leq \beta\ \ps$, and $\as$,
$$\RE\ y> 0,\  \RE\ z\in \R^d,\ \ -\gamma u(t)y-{|z|^2\over y}\leq G(\omega,t,y,z)\leq \gamma u(t)y.$$
Furthermore, for each pair of positive real numbers $K_1$ and $K_2$ satisfying that $[Q_0,S_0]\subset [2K_1,K_2/2]$, let $\Psi$ be a smooth function with values in $[0,1]$ which satisfies that $\Psi(x)=1$ for $x\in [2K_1,K_2/2]$, and $\Psi(x)=0$ when $x$ is outside $[K_1,K_2]$. Define
$$G_\Psi(\omega,t,y,z):=\Psi (y)G(\omega,t,y,z),$$
then $\as$,
$$\RE\ y> 0,\  \RE\ z\in \R^d,\ \ -\gamma u(t)y-\left({1\over K_1}\right)|z|^2\leq G_\Psi(\omega,t,y,z)\leq \gamma u(t)y.$$
Thanks to Lemma 3.3, we know that the BSDE
$$
Y_t=\eta+\int_t^T G_\Psi(s,Y_s,Z_s)\ {\rm d}s-\int_t^T Z_s\cdot {\rm d}B_s,\ \ t\in\T
$$
has a maximal one $(Y^\Psi,Z^\Psi)$ among all bounded solutions. Moreover, we also have that for each $t\in \T$, $0<Q_0\leq Q_t\leq Y_t^{\Psi}\leq S_t\leq S_0\ \ps$, which means that, in view of the definition of $\Psi$, $(Y^\Psi,Z^\Psi)$ is a bounded solution to the following BSDE:
$$
Y_t=\eta+\int_t^T G(s,Y_s,Z_s)\ {\rm d}s-\int_t^T Z_s\cdot {\rm d}B_s,\ \ t\in\T.
$$
We define
$$y^\Psi:={\ln(Y^\Psi)\over \gamma}\ \  {\rm and}\ \ z^\Psi:={Z^\Psi\over \gamma Y^\Psi}.\vspace{0.1cm}$$
It follows from It\^{o}'s formula that $(y^\Psi,z^\Psi)$ is a bounded solution of BSDE$(\xi,g)$, and it is easy to verify that for each $t\in \T$,
\begin{equation}
|y^\Psi_t|\leq \|\xi\|_{\infty}+\int_t^T u(s)\ {\rm d}s\ \ \ \ps.\vspace{-0.2cm}
\end{equation}

In the sequel, let us show that $(y^\Psi,z^\Psi)$ is also the maximal bounded solution. Indeed, let $(\hat y, \hat z)$ be a bounded solution of BSDE$(\xi,g)$, with $A\leq \hat y\leq B$. We can choose positive real numbers $\hat K_1,\hat K_2$ satisfying that $[e^{\gamma A},e^{\gamma B}]\subset [2\hat K_1,\hat K_2]$ and consider $\hat \Psi$ with values in $[0,1]$ which satisfies that $\hat\Psi(x)=1$ for $x\in [2\hat K_1,\hat K_2/2]$, and $\hat \Psi(x)=0$ when $x$ is outside $[\hat K_1,\hat K_2]$. It then follows from It\^{o}'s formula that $\hat Y:=e^{\gamma \hat y},\hat Z:=\gamma \hat z \hat y$ is a bounded solution of BSDE$(\eta,G_{ \hat \Psi})$, where $G_{\hat \Psi}(\omega,t,y,z):=\hat \Psi (y)G(\omega,t,y,z)$. Note that $\as$,
$$\RE\ y> 0,\  \RE\ z\in \R^d,\ \ -\gamma u(t)y-\left({1\over \hat K_1}\right)|z|^2\leq G_{\hat \Psi}(\omega,t,y,z)\leq \gamma u(t)y.$$
Thanks to Lemma 3.3 again, we know that $0<Q_0\leq Q_t\leq \hat Y_t\leq S_t\leq S_0$, which means that $(\hat Y,\hat Z)$ is also a bounded solution of BSDE$(\eta,G_\Psi)$. Therefore, since $(Y^\Psi,Z^\Psi)$ is the maximal bounded solution of BSDE$(\eta,G_\Psi)$, we have $e^{\gamma \hat y}=\hat Y\leq Y^\Psi=e^{\gamma y^\Psi}$, and then for each $t\in \T$, $\hat y_t\leq y^\Psi_t\ \ \ps$. This is the desired result. This argument also shows that $y^\Psi$ does not depend on $\Psi$.

Finally, let us define
$$\tilde g(\omega,t,y,z):=-g(\omega,t,-y,-z),\ \ \RE\ \omega,t,y,z.$$
Then $\tilde g$ also satisfies (2A3) and (25). Consequently, by above arguments we know that BSDE$(-\xi,\tilde g)$ has a maximal bounded solution $(\tilde y,\tilde z)$ and $\tilde y$ also satisfies the estimate in (26). Furthermore, it is easy to verify that $(-\tilde y,-\tilde z)$ is just the minimal bounded solution of BSDE$(\xi,g)$. Lemma 3.4 is then proved.\vspace{0.2cm}\hfill $\Box$

{\bf Remark 3.1} A similar result to Lemma 3.4 was given in \citet{Mor09}, but a different method is used there.\vspace{0.2cm}

We are now in a position to prove Theorem 3.1.\vspace{0.1cm}

{\bf Proof of Theorem 3.1.}\ Assume that $\xi\in L^\infty(\Omega,\F_T,P)$ and $g$ satisfies (2A3), (3A1) and (3A2). We only prove the case of the maximal solution, another case can be proved in a similar way. First, by Lemma 3.1 we can let $(L_t)_{t\in \T}$ and $(U_t)_{t\in\T}$ be, respectively, the unique global solution to the following two backward ODEs
\begin{equation}
L_t=-\|\xi\|_{\infty}-\int_t^T u(s)\ l(L_s)\ {\rm d}s
\end{equation}
and
\begin{equation}
U_t=\|\xi\|_{\infty}+\int_t^T u(s)\ l(U_s)\ {\rm d}s
\end{equation}
Then for each $t\in \T$, we have
$$L_0\leq L_t\leq -\|\xi\|_{\infty}\leq 0\leq \|\xi\|_{\infty}\leq U_t\leq U_0.\vspace{-0.1cm}$$

For each constant $K>0$ satisfying that $[L_0,U_0]\subset [-K,K]$, consider a continuous function $\kappa$ such that $\kappa(x)=-K$ when $x<-K$, $\kappa(x)=x$ when $x\in [-K,K]$, and $\kappa(x)=K$ when $x> K$. For each $(\omega,t,y,z)\in \Omega\times \T\times \R\times \R^d$, we define
$$g_\kappa(\omega,t,y,z):=g(\omega,t,\kappa(y),z)\ \ {\rm and}\ \ \gamma^K:=2\left(\max\limits_{x\in [-K,K]}h(x)+1\right).$$
Then it follows from (3A1) and (3A2) that $\as$, for each $(y,z)\in \R\times \R^d$,
\begin{equation}
g_\kappa(t,y,z)\ {\rm sgn}(y)\leq u(t)l(\kappa(y))+{\gamma^K \over 2}|z|^2
\end{equation}
and
$$
|g_\kappa(t,y,z)|\leq \bar u(t)\left(\max\limits_{x\in [-K,K]}\bar\varphi(x)\right)+\left(\max\limits_{x\in [-K,K]}\bar h(x)\right)|z|^2.\vspace{0.2cm}
$$
It then follows from Lemma 3.4 that BSDE$(\xi,g_\kappa)$ has a maximal bounded solution $(y^\kappa_t,z^\kappa_t)_{t\in \T}$. Let $Y^\kappa_t:=e^{\gamma^K y^\kappa_t}$ and $Z^\kappa_t:=\gamma^K Y^\kappa_tz^\kappa_t$, then $(Y^\kappa,Z^\kappa)$ is a bounded solution to the following BSDE
$$
Y_t=\eta+\int_t^T G_\kappa(s,Y_s,Z_s)\ {\rm d}s-\int_t^T Z_s\cdot {\rm d}B_s,\ \ t\in\T,\vspace{0.1cm}
$$
where $\eta:=e^{\gamma^K \xi}$ and for each $(\omega,t,y,z)\in \Omega\times\T\times\R\times\R^d$,
$$
G_\kappa(\omega,t,y,z):=\Dis \left[\gamma^K yg_\kappa\left(\omega,t,{\ln y\over \gamma^K},{z\over \gamma^K y}\right)-{1\over 2}{|z|^2\over y}\right]\mathbbm{1}_{ y>0}.
$$
Furthermore, it follows from (29) that $\as$,
\begin{equation}
\RE\ y>0,\ \ \RE\ z\in\R^d,\ \ g_\kappa(t,y,z)\leq u(t)l(\kappa(y))+{\gamma^K \over 2}|z|^2,
\end{equation}
and then $\as$,
\begin{equation}
\RE\ y>1,\ \RE\ z\in \R^d,\ \ G_\kappa(\omega,t,y,z)\leq G'_\kappa(\omega,t,y,z):=\gamma^K yu(t)l\left(\kappa\left({\ln y\over \gamma^K}\right)\right)\mathbbm{1}_{ y>0}.
\end{equation}
Note that $\as$, for each $(y,z)\in \R\times\R^d$,
\begin{equation}
|G'_\kappa(t,y,z)|\leq \gamma^K \left(\max\limits_{x\in [-K,K]}l(x)\right) u(t)|y|.
\end{equation}
By Lemma 2.2 we know that BSDE$(\|\eta\|_\infty,G'_\kappa)$ has a maximal $L^2$ solution $({}^\kappa Y',{}^\kappa Z')$. On the other hand, noticing the definition of $\kappa$ and the assumptions of $u(\cdot)$ and $l(\cdot)$, we can verify directly that the following backward ODE
$$R^K_t=e^{\gamma^K\|\xi\|_{\infty}}+\int_t^T \gamma^K u(s)l\left(\kappa\left({\ln R^K_s\over \gamma^K}\right)\right)R^K_s\mathbbm{1}_{ R^K_s>0}\ {\rm d}s$$
has a unique solution $R^K_t=e^{\gamma^K U_t}$ with $t\in \T$, where $U_t$ is defined in (28). Thus, by terms of Lemma 3.2 we know that $({}^\kappa Y',{}^\kappa Z')=(R^K,0)$. Furthermore, in view of (31), (32) and the fact ${}^\kappa Y'\geq 1$, it follows from Theorem 2.2 and Remark 2.2 that for each $t\in \T$,
$$e^{\gamma^K y^\kappa_t}=Y^\kappa_t\leq {}^\kappa Y'_t=e^{\gamma^K U_t}\ \ \ps,$$
and then
\begin{equation}
y^\kappa_t\leq U_t\leq U_0\leq K\ \ \ps.
\end{equation}

On the other hand, noticing that $\kappa(-x)=-\kappa(x)$, by (29) we get that $\as$,
$$\RE\ y> 0,\ \ \RE\ z\in\R^d,\ \ \bar g_\kappa(\omega,t,y,z):=-g_\kappa(\omega,t,-y,-z)\leq u(t)\bar l(\kappa(y))+{\gamma^K \over 2}|z|^2,$$
where $\bar l(u):=l(-u)$ for each $u\in \R$. Furthermore, note that $(-y^\kappa_t,-z^\kappa_t)_{t\in \T}$ is a bounded solution of BSDE$(-\xi,\bar g_\kappa)$, $\bar l(\cdot)\in \LS$ and the following backward ODE
$$
\bar U_t=\|-\xi\|_{\infty}+\int_t^T u(s)\ \bar l(\bar U_s)\ {\rm d}s=\|\xi\|_{\infty}+\int_t^T u(s)\ l(-\bar U_s)\ {\rm d}s
$$
has a unique solution $\bar U_t=-L_t$ with $t\in \T$, where $L_t$ is defined in (27). The similar argument to that from (29) to (33) yields that for each $t\in \T$,
$$
-y^\kappa_t\leq \bar U_t=-L_t\ \ \ps,
$$
and then
\begin{equation}
y^\kappa_t\geq L_t\geq L_0\geq -K\ \ \ps.\vspace{0.2cm}
\end{equation}
Thus, in view of (33), (34), the definition of $g_\kappa$ and $\kappa$, we know that the $(y^\kappa_t,z^\kappa_t)_{t\in \T}$ is a bounded solution of BSDE$(\xi,g)$.

Finally, let us show that $(y^\kappa_t,z^\kappa_t)_{t\in \T}$ is also the maximal bounded solution of BSDE$(\xi,g)$. Indeed, let $(\hat y_t, \hat z_t)_{t\in \T}$ be a bounded solution of BSDE$(\xi,g)$, with $A\leq \hat y\leq B$. We choose a positive constant $\hat K$ satisfying that $[A,B]\cup[L_0,U_0]\subset [-\hat K,\hat K]$ and consider $\hat \kappa$ such that $\hat\kappa(x)=-\hat K$ for $x<-\hat K$, $\hat\kappa(x)=x$ for $x\in [-\hat K,\hat K]$, and $\hat\kappa(x)=\hat K$ for $x>\hat K$. For each $(\omega,t,y,z)\in \Omega\times \T\times \R\times \R^d$, we define
$$g_{\hat\kappa}(\omega,t,y,z):=g(\omega,t,{\hat\kappa(y)},z).$$
Then $(\hat y_t, \hat z_t)_{t\in \T}$ is a bounded solution of BSDE$(\xi,g_{\hat\kappa})$. Furthermore, by the above argument as above (from (29) to (34)) we can conclude that for each $t\in \T$,
$$-K\leq L_0\leq L_t\leq \hat y_t\leq U_t\leq U_0\leq K\ \ \ps,$$
which means that $(\hat y_t, \hat z_t)_{t\in \T}$ is also a bounded solution of BSDE$(\xi,g_{\kappa})$. Note that $(y^\kappa_t, z^\kappa_t)_{t\in \T}$ is the maximal bounded solution of BSDE$(\xi,g_{\kappa})$. We know that for each $t\in \T$, $\hat y_t\leq y^\kappa_t\ \ps$, which is the desired result. This argument also shows that $y^\kappa_t$ does not depend on $\kappa$. The proof of Theorem 3.1 is then complete.\vspace{0.2cm}\hfill $\Box$

Finally, in view of Theorems 2.2 and 2.3, by checking carefully the proof of Theorem 3.1, Lemma 3.4 and, especially, Lemma 3.3, we can prove that the following two comparison theorems for the maximal and minimal bounded solutions hold true. \vspace{0.1cm}

{\bf Theorem 3.2}\ \ Assume that $0<T\leq +\infty$, $g$ and $g'$ are two generators of BSDEs, and $(y_\cdot,z_\cdot)$ is any bounded solution of BSDE$(\xi,g)$. Assume further that $g'$ satisfies (2A3), (3A1) and (3A2), and $(y'_\cdot,z'_\cdot)$ is the maximal bounded solution of BSDE$(\xi',g')$ by Theorem 3.1. If $\xi\leq \xi'\ \ps$ and (9) holds true, then for each $t\in\T$,
$$y_t\leq y'_t\ \ \ \ps.$$

{\bf Theorem 3.3}\ Assume that $0<T\leq +\infty$, $g$ and $g'$ are two generators of BSDEs, and $(y'_\cdot,z'_\cdot)$ is any bounded solution of BSDE$(\xi',g')$. Assume further that $g$ satisfies (2A3), (3A1) and (3A2), and $(y_\cdot,z_\cdot)$ is the minimal bounded solution of BSDE$(\xi,g)$ by Theorem 3.1. If $\xi\leq \xi'\ \ps$ and (8) holds true, then for each $t\in\T$,
$$y_t\leq y'_t\ \ \ \ps.$$

{\bf Corollary 3.1} Assume that $0<T\leq +\infty$ and both $g$ and $g'$ satisfy (2A3), (3A1) and (3A2). Let $(y_\cdot,z_\cdot)$ and $(y'_\cdot,z'_\cdot)$ be, respectively, the maximal (resp. minimal) bounded solution of BSDE$(\xi,g)$ and BSDE$(\xi',g')$ by Theorem 3.1. If $\xi\leq \xi'\ \ps$ and for each $(y,z)\in \R\times \R^d$,
$$g(t,y,z)\leq g'(t,y,z)\ \ \ \as,$$
then for each $t\in \T$, we have
$$y_t\leq y'_t\ \ \ \ps.$$

\section{Comparison theorems of bounded solutions}

In this section, we will establish two comparison theorems and an existence and uniqueness theorem for bounded solutions of BSDEs. Let us first introduce the following assumptions on the generator $g$, where $0<T\leq +\infty$:\vspace{0.2cm}

{\bf (4A1)} There exists a function $v(\cdot)\in \VS$ such that $\as$,
$$|g(\omega,t,y,z_1)-g(\omega,t,y,z_2)|\leq (v(t)+|z_1|+|z_2|)|z_1-z_2|,\ \ \RE\ y,z_1,z_2.$$

{\bf (4A2)} $\as$, $\RE\ y$, $g(\omega,t,y,\cdot):\R^d\mapsto \R$ is convex, or $\as$, $\RE\ y$, $g(\omega,t,y,\cdot):\R^d\mapsto \R$ is concave. \vspace{0.2cm}

Let us now recall several facts on the martingales of bounded mean oscillation, briefly called BMO-martingales. Readers are refereed to \citet{Kaz94}, \citet{Hu05} and \citet{BC08} for more details. Suppose that $(z_t)_{t\in\T}$ is an $(\F_t)$-progressively measurable $\R^d$-valued process such that $\int_0^T |z_s|^2\ {\rm d}s<+\infty$, $\ps$. First, it is well known that the local martingale $\int_0^\cdot z_s\cdot {\rm d}B_s$ is a BMO-martingale if and only if
$$\sup\limits_{\tau \in\Sigma_T} \E\left[\int_{\tau}^T |z_s|^2\ {\rm d}s|\F_\tau\right]<+\infty\ \ \ \ps.\vspace{0.1cm}$$
Furthermore, according to Theorem 2.3 in \citet{Kaz94}, the stochastic exponential ${\mathcal E}(M)$ of a BMO-martingale $M$ is a uniformly integrable martingale, note that the stochastic exponential ${\mathcal E}(M)$ is given by
$${\mathcal E}(M)_t=\exp(M_t-{1\over 2}\langle M\rangle_t), \ \ t\in\T,\vspace{0.1cm}$$
where the quadratic variation is denoted by $\langle M\rangle$. Finally, by Theorem 3.6 in \citet{Kaz94} we also know that if $Q$ is a probability measure defined by
$${\rm d}Q={\mathcal E}(M)_T{\rm d}P$$
for a BMO-martingale $M$ under $P$, then the Girsanov transform of a BMO-martingale under $P$ is a BMO-martingale under $Q$.\vspace{0.1cm}

The following Lemma 4.1 will be used several times later.\vspace{0.1cm}

{\bf Lemma 4.1}\ Let $0<T\leq +\infty$, $g$ satisfies (3A2), and let $(y_\cdot,z_\cdot)$ be a bounded solution of BSDE$(\xi,g)$. Then $\int_0^\cdot z_s\cdot {\rm d}B_s$ is a BMO-martingale under $P$.\vspace{0.1cm}

{\bf Proof.}\ Assume that $A\leq y_\cdot\leq B$ with $A,B\in \R$. Let $\gamma=2\left(\max\limits_{x\in [A,B]}\bar h(x)+1\right)$ and consider the following function from $\R_+$ into itself defined by
$$f(x)={1\over \gamma^2}(e^{\gamma x}-1-\gamma x).$$
Clearly, $x\mapsto f(|x|)$ is ${\mathcal C}^2$ and for each $\tau\in \Sigma_T$, we have from It\^{o}'s formula,
$$\begin{array}{lll}
f(|y_\tau|)&=& \Dis f(|\xi|)+\int_\tau^T \left(f'(|y_s|){\rm sgn}(y_s)g(s,y_s,z_s)-{1\over 2}f''(|y_s|)|z_s|^2\right){\rm d}s\\
&&\Dis -\int_\tau^T f'(|y_s|){\rm sgn}(y_s)z_s\cdot {\rm d}B_s.
\end{array}$$
Since $A\leq y_\cdot\leq B$ and $f'(x)\geq 0$ for $x\geq 0$, by (3A2) and the definition of $\gamma$ we can get the existence of a constant $k>0$ satisfying that
$$\begin{array}{lll}
0\leq f(|y_\tau|)&\leq & \Dis k\left(1+\int_0^T \bar u(s)\ {\rm d}s\right)-\int_\tau^T f'(|y_s|){\rm sgn}(y_s)z_s\cdot {\rm d}B_s\\
&&\Dis -{1\over 2}\int_\tau^T \left[\left(f''(|y_s|)-\gamma f'(|y_s|)\right)|z_s|^2\right]\ {\rm d}s.
\end{array}$$
Note that $(y_\cdot,z_\cdot)$ be a bounded solution and $f''(x)-\gamma f'(x)=1$ for $x\geq 0$. By taking the conditional expectation with respect to $\F_\tau$ under $P$ in the previous inequality we get that for each $\tau\in \Sigma_T$,
$$\E\left[\left.\int_\tau^T |z_s|^2\ {\rm d}s\right|\F_\tau\right]\leq 2k\left(1+\int_0^T \bar u(s)\ {\rm d}s\right).$$
That is to say, $\int_0^\cdot z_s\cdot {\rm d}B_s$ is a BMO-martingale under $P$. The proof is complete. \vspace{0.2cm}\hfill $\Box$

The following Theorems 4.1-4.2 establish two comparison theorems for bounded solutions of BSDEs, which virtually improves the corresponding comparison results obtained in \citet{Bri08} and \citet{Mor09} even for the case of the finite time horizon.\vspace{0.1cm}

{\bf Theorem 4.1} \ Assume that $0<T\leq +\infty$, both $g$ and $g'$ satisfy (3A2), and $(y_\cdot,z_\cdot)$ and $(y'_\cdot,z'_\cdot)$ are, respectively, a bounded solution of BSDE$(\xi,g)$ and BSDE$(\xi',g')$. If $\xi\leq \xi'\ \ps$ and one of the following two statements is satisfied:\vspace{0.2cm}

(i)\ \ $g$ satisfies (2A1) and (4A1), and (8) holds true;\vspace{0.1cm}

(ii) $g'$ satisfies (2A1) and (4A1), and (9) holds true,\vspace{0.2cm}

\noindent then for each $t\in\T$,
$$y_t\leq y'_t\ \ \ \ps.$$

{\bf Proof.}\ We only prove the case (i). Another case can be proved in the same way. Tanaka's formula leads to the equation, setting
$\hat{y}_t=y_t-y'_t,\ \hat{z}_t=z_t-z'_t$ and noticing that $\ps,\ (\xi-\xi')^+=0$,
\begin{equation}
\Dis \hat{y}_t^+\leq  \int_t^T
 \mathbbm{1}_{ \hat{y}_s> 0}(g(s,y_s,z_s)-g'(s,y'_s,z'_s))
 \ {\rm d}s-\int_t^T\mathbbm{1}_{ \hat{y}_s> 0}
\hat{z}_s\cdot {\rm d}B_s,\ \ t\in \T.
\end{equation}
First of all, since $\mathbbm{1}_{ \hat{y}_s> 0}(g(s,y'_s,z'_s)-g'(s,y'_s,z'_s))$ is
non-positive, we have
$$
\begin{array}{lll}
&&\mathbbm{1}_{ \hat{y}_s> 0}(g(s,y_s,z_s)-g'(s,y'_s,z'_s))\\
&=&\mathbbm{1}_{ \hat{y}_s> 0}(g(s,y_s,z_s)-g(s,y'_s,z'_s))
+\mathbbm{1}_{ \hat{y}_s> 0}(g(s,y'_s,z'_s)-g'(s,y'_s,z'_s))\\
&\leq & \mathbbm{1}_{ \hat{y}_s> 0}(g(s,y_s,z_s)-g(s,y'_s,z_s))+\mathbbm{1}_{ \hat{y}_s> 0}(g(s,y'_s,z_s)
-g(s,y'_s,z'_s))
\end{array}$$
and we deduce, using assumptions (2A1) and (4A1) for $g$, that
\begin{equation}
\mathbbm{1}_{ \hat{y}_s> 0}(g(s,y_s,z_s)-g'(s,y'_s,z'_s))\leq
u(s)\rho(\hat{y}_s^+)+\mathbbm{1}_{ \hat{y}_s> 0}(v(s)+|z_s|+|z'_s|)|\hat{z}_s|.
\end{equation}
Thus, by (35) and (36) we can get that for each $t\in\T$,
\begin{equation}
\begin{array}{lll}
\Dis \hat{y}_t^+ &\leq & \Dis\int_t^T
\left[u(s)\rho(\hat{y}_s^+)+\mathbbm{1}_{ \hat{y}_s>
0}(v(s)+|z_s|+|z'_s|)|\hat{z}_s|\right]{\rm d}s -\int_t^T \mathbbm{1}_{ \hat{y}_s> 0}\hat{z}_s\cdot {\rm d}B_s\\
&=&\Dis \int_t^Tu(s)\rho(\hat{y}_s^+)\ {\rm
d}s-\int_t^T \mathbbm{1}_{ \hat{y}_s> 0}\hat{z}_s\cdot [{\rm d}B_s-{(v(s)+|z_s|+|z'_s|)\hat{z}_s\over
|\hat{z}_s|}\mathbbm{1}_{ |\hat{z}_s|\neq 0}\ {\rm d}s].
\end{array}
\end{equation}
Furthermore, since  both $g$ and $g'$ satisfy (3A2), and both $(y_\cdot,z_\cdot)$ and $(y'_\cdot,z'_\cdot)$ are bounded solutions, it follows from Lemma 4.1 that both $\int_0^\cdot z_s\cdot {\rm d}B_s$ and $\int_0^\cdot z'_s\cdot {\rm d}B_s$ are BMO-martingale under $P$. Then we have
$$
\begin{array}{lll}
&&\Dis \sup\limits_{\tau\in \Sigma_T}\E\left[\int_\tau^T (v(s)+|z_s|+|z'_s|)^2\ {\rm d}s|\F_\tau\right]\\
&\leq &\Dis 4\int_0^T v^2(s)
\ {\rm d}s+4\sup\limits_{\tau\in \Sigma_T}\E\left[\int_\tau^T |z_s|^2\ {\rm d}s|\F_\tau\right]+4\sup\limits_{\tau\in \Sigma_T}\E\left[\int_\tau^T |z'_s|^2\ {\rm d}s|\F_\tau\right]<+\infty,
\end{array}
$$
which means that the process
$$M_t:=\int_0^t {(v(s)+|z_s|+|z'_s|)\hat{z}_s\over
|\hat{z}_s|}\mathbbm{1}_{ |\hat{z}_s|\neq 0}\ \cdot {\rm d}B_s,\ \ t\in \T$$
is a BMO-martingale under $P$. Then the stochastic exponential $${\mathcal E}(M)_t=\exp(M_t-{1\over 2}\langle M\rangle_t), \ \ t\in\T$$
of $M$ is a uniformly integrable martingale. Now, let us define by $Q$ the probability measure under $(\Omega,\F_T)$ given by
$${{\rm d}Q\over {\rm d}P}:={\mathcal E}(M)_T.$$
Then, noticing that
$$\overline {M}_t:=\int_0^t \mathbbm{1}_{ \hat{y}_s> 0}\hat{z}_s\cdot {\rm d}B_s,\ \ t\in \T\vspace{0.1cm}$$
is a also BMO-martingale under $P$, we know that the process
$$\int_0^t \mathbbm{1}_{ \hat{y}_s> 0}\hat{z}_s\cdot [{\rm d}B_s-{(v(s)+|z_s|+|z'_s|)\hat{z}_s\over
|\hat{z}_s|}\mathbbm{1}_{ |\hat{z}_s|\neq 0}\ {\rm d}s],\ \ t\in \T,$$
the Girsanov transform of $\overline {M}$, is a BMO-martingale under $Q$. Let $\E_Q[X|\F_t]$ represent the conditional expectation of the random variable $X$ with respect to $\F_t$ under $Q$. Taking the conditional expectation with respect to $\F_t$ under $Q$ in (37) yields that for each $t\in\T$,
$$\hat y_t^+\leq \E_Q\left[\left.\int_t^T u(s) \rho(\hat y_s^+)\ {\rm d}s\right|\F_t\right]\ \ \ps.\vspace{0.1cm}$$
Thus, applying Lemma 2.1 with $u_n(t)\equiv \hat{y}_t^+$,
$b_n\equiv 0$, $P_n\equiv Q$, $\beta (s)=u(s)$ and $\psi(u)=\rho(u)$ yields that for each $t\in \T$,
$$\hat{y}_t^+=\Lim \hat{y}_t^+=\Lim u_n(t)=0\ \ \ \ps.$$
That is to say, for each $t\in\T$, $\ y_t\leq y'_t\ \ \ps$. The proof is complete.\vspace{0.2cm}\hfill $\Box$

{\bf Theorem 4.2} \ Assume that $0<T\leq +\infty$, both $g$ and $g'$ satisfy (3A2), and $(y_\cdot,z_\cdot)$ and $(y'_\cdot,z'_\cdot)$ are, respectively, a bounded solution of BSDE$(\xi,g)$ and BSDE$(\xi',g')$. If $\xi\leq \xi'\ \ps$ and one of the following two statements is satisfied:\vspace{0.2cm}

(i)\ \ $g$ satisfies (2A1) and (4A2), and (8) holds true;\vspace{0.1cm}

(ii) $g'$ satisfies (2A1) and (4A2), and (9) holds true,\vspace{0.2cm}

\noindent then for each $t\in\T$,
$$y_t\leq y'_t\ \ \ \ps.$$

{\bf Proof.}\ The proof will be split into four steps.

{\bf First step:}\ \ Suppose that both $g$ and $g'$ satisfy (3A2), $\xi\leq \xi'\ \ps$, $g$ also satisfies (2A1) and is convex with respect to $z$, and (8) holds true. We further assume that both $y_\cdot$ and $y'_\cdot$ is non-positive.

For each $n\geq 2$, let us set
$$\hat{y}^n_t=y_t-{n-1\over n} y'_t,\ \hat{z}^n_t=z_t-{n-1\over n} z'_t.$$
Note that
$$(\xi-{n-1\over n}\xi')^+\leq (\xi-\xi')^+=0$$
due to the fact that $\xi'\leq 0$. Tanaka's formula yields that for $t\in \T$,
\begin{equation}
\Dis (\hat{y}^n_t)^+\leq  \int_t^T
 \mathbbm{1}_{ \hat{y}^n_s> 0}\left[g(s,y_s,z_s)-{n-1\over n}g'(s,y'_s,z'_s)\right]
 \ {\rm d}s-\int_t^T\mathbbm{1}_{ \hat{y}^n_s> 0}
\hat{z}^n_s\cdot {\rm d}B_s.
\end{equation}
First of all, in view of $y'_s \leq 0$ and then $${n-1\over n}y'_s\geq y'_s,$$
it follows from (8) that $\mathbbm{1}_{ \hat{y}^n_s> 0}(g(s,y'_s,z'_s)-g'(s,y'_s,z'_s))$ is
non-positive. Then we have
\begin{equation}
\begin{array}{lll}
&& \Dis \mathbbm{1}_{ \hat{y}^n_s> 0}(g(s,y_s,z_s)-{n-1\over n}g'(s,y'_s,z'_s))\vspace{0.1cm}\\
&=&\Dis \mathbbm{1}_{ \hat{y}^n_s> 0}(g(s,y_s,z_s)-{n-1\over n}g(s,y'_s,z'_s))
+\mathbbm{1}_{ \hat{y}^n_s> 0}{n-1\over n}(g(s,y'_s,z'_s)-g'(s,y'_s,z'_s))\vspace{0.1cm}\\
&\leq & \Dis \mathbbm{1}_{ \hat{y}^n_s> 0}(g(s,y_s,z_s)-g(s,y'_s,z_s))+\mathbbm{1}_{ \hat{y}^n_s> 0}(g(s,y'_s,z_s)
-{n-1\over n}g(s,y'_s,z'_s)).
\end{array}
\end{equation}
By (2A1) we can obtain that, in view of $y'_s\leq 0$,
\begin{equation}
\begin{array}{lll}
&&\Dis \mathbbm{1}_{ \hat{y}^n_s> 0}(g(s,y_s,z_s)-g(s,y'_s,z_s))\vspace{0.2cm}\\
&= &\Dis \mathbbm{1}_{ \hat{y}^n_s> 0}\left[(g(s,y_s,z_s)-g(s,{n-1\over n}y'_s,z_s))+(g(s,{n-1\over n}y'_s,z_s)-g(s,y'_s,z_s))\right]\\
&\leq & \Dis u(s)\rho((\hat{y}^n_s)^+)+u(s)\rho({-y'_s\over n}).
\end{array}
\end{equation}
Furthermore, since $g$ is convex with respect to $z$, and satisfies (3A2), we have
\begin{equation}
\begin{array}{lll}
g(s,y'_s,z_s)&=& \Dis g\left(s, y'_s,{n-1\over n}z'_s+{1\over n}( nz_s-(n-1)z'_s)\right)\vspace{0.2cm}\\
&\leq & \Dis {n-1\over n} g(s,y'_s,z'_s)+{1\over n}g\left(s,y'_s,n\hat z^n_s\right)\vspace{0.2cm}\\
&\leq & \Dis {n-1\over n} g(s,y'_s,z'_s)+{\bar u(s)\bar\varphi(y'_s)\over n}+n\bar h(y'_s)|\hat z^n_s|^2.
\end{array}
\end{equation}
Thus, in view of the fact that $-k\leq y'_s\leq 0$ for some positive constant $k>0$, combining (38)-(41) yields that
\begin{equation}
\begin{array}{lll}
\Dis (\hat{y}^n_t)^+ &\leq & \Dis a_n +\int_t^T
 (u(s)\rho((\hat{y}^n_s)^+)+n\gamma\mathbbm{1}_{ \hat{y}^n_s> 0}|\hat{z}_t^n|^2) \ {\rm d}s-\int_t^T \mathbbm{1}_{ \hat{y}^n_s> 0}
\hat{z}^n_s\cdot {\rm d}B_s\\
&=& \Dis a_n +\int_t^T
 u(s)\rho((\hat{y}^n_s)^+)\ {\rm d}s-\int_t^T{\mathbbm{1}_{ \hat{y}^n_s> 0}
\hat{z}^n_s\cdot [\rm d}B_s-n\gamma \hat{z}_t^n\ {\rm d}s],\ t\in\T,
\end{array}
\end{equation}
where $\gamma:=\max\limits_{x\in [-k,0]}\bar h(x)$ and
$$a_n:=\rho({k\over n})\int_0^T u(s)\ {\rm d}s+{\max\limits_{x\in [-k,0]}\bar\varphi(x)\over n}\int_0^T \bar u(s)\ {\rm d}s\To 0\ \ {\rm as}\ \ n\To \infty.\vspace{0.1cm}$$

In the sequel, since both $g$ and $g'$ satisfy (3A2), and both $(y_\cdot,z_\cdot)$ and $(y'_\cdot,z'_\cdot)$ are bounded solutions, it follows from Lemma 4.1 that for each $n\geq 2$, the process
$$N_t^n:=\int_0^t n\gamma \hat{z}_t^n\ \cdot {\rm d}B_s,\ \ t\in \T\vspace{0.1cm}$$
is a BMO-martingale under $P$. Then the stochastic exponential $${\mathcal E}(N^n)_t=\exp(N_t^n-{1\over 2}\langle N^n\rangle_t), \ \ t\in\T$$
of $N^n$ is a uniformly integrable martingale. Now, let us define by $P_n$ the probability measure under $(\Omega,\F_T)$ given by
$${{\rm d}P_n\over {\rm d}P}:={\mathcal E}(N^n)_T.$$
Then, noticing that
$$\overline{N}^n_t:=\int_0^t \mathbbm{1}_{ \hat{y}^n_s> 0}\hat{z}^n_s\cdot {\rm d}B_s,\ \ t\in \T\vspace{0.1cm}$$
is a also BMO-martingale under $P$, we know that the process
$$\int_0^t \mathbbm{1}_{ \hat{y}^n_s> 0}\hat{z}^n_s\cdot [{\rm d}B_s-n\gamma \hat{z}_t^n\ {\rm d}s],\ \ t\in \T,$$
the Girsanov transform of $\overline{N}^n$, is a BMO-martingale under $P_n$. Let $\E_n[X|\F_t]$ represent the conditional expectation of the random variable $X$ with respect to $\F_t$ under $P_n$. Taking the conditional expectation with respect to $\F_t$ under $P_n$ in (42) yields that for each $t\in\T$,
$$(\hat{y}_t^n)^+\leq a_n+\E_n\left[\left.\int_t^T u(s) \rho((\hat{y}_s^n)^+)\ {\rm d}s\right|\F_t\right]\ \ \ps.$$
Thus, applying Lemma 2.1 with $u_n(t)=(\hat{y}_t^n)^+$, $b_n=a_n$, $\beta (s)=u(s)$ and $\psi(u)=\rho(u)$ yields that for each $t\in\T$,
$$(y_t-y'_t)^+=\Lim (\hat{y}_t^n)^+=\Lim u_n(t)=0\ \ \ \ps.$$
That is to say, for each $t\in\T$, $\ y_t\leq y'_t\ \ \ps$.\vspace{0.2cm}

{\bf Second step:}\ In this step, we will eliminate the condition that both $y_\cdot$ and $y'_\cdot$ is non-positive required in the first step.

Indeed, note that both $y_\cdot$ and $y'_\cdot$ are bounded processes. We can assume that there exists a constant $k>0$ such that for each $t\in \T$, $|y_t|+|y'_t|\leq k\ \ps$. Let $\bar y_t:=y_t-k$ and $\bar y'_t:=y'_t-k$, then $(\bar y_t,z_t)_{t\in\T}$ and $(\bar y'_t,z'_t)_{t\in\T}$ are, respectively, a bounded solution of BSDE$(\bar\xi,\bar g)$ and BSDE$(\bar\xi',\bar g')$, where $\bar\xi:=\xi-k$, $\bar\xi':=\xi'-k$ and for each $(\omega,t,y,z)\in \Omega\times\T\times\R\times \R^d$,
$$\bar g(\omega,t,y,z):=g(\omega,t,y+k,z),\ \ \bar g'(\omega,t,y,z):=g'(\omega,t,y+k,z).$$
It is not difficult to verify that both $\bar g$ and $\bar g'$ satisfy (3A2) but $\bar\varphi(u)$ and $\bar h(u)$ are replaced by
$\bar\varphi(u+k)$ and $\bar h(u+k)$ respectively, $\bar\xi\leq \bar\xi'\ \ps$, $\bar g$ satisfies (2A1), $\as$,
$$\begin{array}{lll}
\Dis \mathbbm{1}_{ \bar y_t>\bar y'_t}\bar g(t,\bar y'_t,z'_t)&=& \Dis \mathbbm{1}_{ y_t>y'_t}g(t,y'_t,z'_t)\\
&\leq & \Dis \mathbbm{1}_{ y_t>y'_t}g'(t,y'_t,z'_t)\\
&=& \Dis \mathbbm{1}_{ \bar y_t>\bar y'_t}\bar g'(t,\bar y'_t,z'_t),
\end{array}$$
$\bar g$ is convex with respect to $z$, and both $\bar y_\cdot$ and $\bar y'_\cdot$ are non-positive. Thus, by the first step we can conclude that for each $t\in \T$,
$$y_t-k=\bar y_t\leq \bar y'_t=y'_t-k\ \ \ \ps,$$
which is the desired result.\vspace{0.1cm}

{\bf Third step:}\ \ Suppose that both $g$ and $g'$ satisfy (3A2), $\xi\leq \xi'\ \ps$, $g'$ also satisfies (2A1) and is convex with respect to $z$, and (9) holds true. Then we can replace (39) by the following inequality
$$
\begin{array}{lll}
&& \mathbbm{1}_{ \hat{y}^n_s> 0}(g(s,y_s,z_s)-{n-1\over n}g'(s,y'_s,z'_s))\vspace{0.1cm}\\
&=&\mathbbm{1}_{ \hat{y}^n_s> 0}(g(s,y_s,z_s)-g'(s,y_s,z_s))
+\mathbbm{1}_{ \hat{y}^n_s> 0}(g'(s,y_s,z_s)-{n-1\over n}g'(s,y'_s,z'_s))\vspace{0.1cm}\\
&\leq & \mathbbm{1}_{ \hat{y}^n_s> 0}(g'(s,y_s,z_s)-g'(s,y'_s,z_s))+\mathbbm{1}_{ \hat{y}^n_s> 0}(g'(s,y'_s,z_s)
-{n-1\over n}g'(s,y'_s,z'_s)),
\end{array}
$$
and then make use of (9) and the assumptions on $g'$ to conclude that (42) holds still true provided that $y'_\cdot\leq 0$. Thus, the similar argument to the previous two steps will give the desired conclusion.\vspace{0.2cm}

{\bf Fourth step:}\  Suppose that both $g$ and $g'$ satisfy (3A2), $\xi\leq \xi'\ \ps$, $g$ also satisfies (2A1) and is concave with respect to $z$, and (8) holds true.

Let us set
$$\tilde \xi':=-\xi',\ \tilde y'_t:=-y'_t,\ \tilde z'_t:=-z'_t,\
\tilde g'(t,y,z):=-g'(t,-y,-z)$$
and
$$\tilde \xi:=-\xi,\ \tilde y_t:=-y_t,\ \tilde z_t:=-z_t,\ \tilde g(t,y,z):=-g(t,-y,-z).\vspace{0.1cm}$$
Then $(\tilde y'_t, \tilde z'_t)_{t\in\T}$ and $(\tilde y_t, \tilde z_t)_{t\in\T}$ are, respectively, a bounded solution of BSDE$(\tilde\xi',\tilde g')$ and BSDE$(\tilde\xi,\tilde g)$. And, it is easy to verify that both $\tilde g'$ and $\tilde g$ satisfy (3A2), $\tilde \xi'\leq \tilde \xi\ \ps$, $\tilde g$ satisfies (2A1) and is convex with respect to $z$, and $\as$,
$$\begin{array}{lll}
\Dis \mathbbm{1}_{ \tilde y'_t>\tilde y_t}\tilde g'(t,\tilde y'_t,\tilde z'_t)&=& \Dis -\mathbbm{1}_{ y_t>y'_t}g'(t,y'_t,z'_t)\\
&\leq & \Dis -\mathbbm{1}_{ y_t>y'_t}g(t,y'_t,z'_t)\\
&= & \Dis \mathbbm{1}_{ \tilde y'_t>\tilde y_t}\tilde g(t,\tilde y'_t,\tilde z'_t).
\end{array}$$
Thus, by the third step we know that for each $t\in \T$,
$$-y'_t=\tilde y'_t\leq \tilde y_t=-y_t\ \ \ps,$$
which is the desired result. In the same way, we can prove the remainder case that both $g$ and $g'$ satisfy (3A2), $\xi\leq \xi'\ \ps$, $g'$ also satisfies (2A1) and is concave with respect to $z$, and (9) holds true. Theorem 4.2 is then proved.\vspace{0.2cm}\hfill $\Box$


By virtue of Theorems 3.1, 2.1, 4.1 and 4.2, we can obtain the following existence and uniqueness result for bounded solutions of BSDEs.\vspace{0.1cm}

{\bf Theorem 4.3}\ Assume that $0<T\leq +\infty$ and that $g$ satisfies (2A3), (3A2) and (2A1). Furthermore, we also assume that one of three assumptions (2A2), (4A1) and (4A2) holds true for $g$. Then for each $\xi\in L^\infty(\Omega,\F_T,P)$, BSDE$(\xi,g)$ has a unique bounded solution $(Y,Z)$. Moreover, for each $t\in \T$, we have
$$L_0\leq L_t\leq Y_t\leq U_t\leq U_0\ \ \ps,$$
where $(L,U)$ are the unique solutions of (27) and (28) but $u(t)$ is replaced by $u(t)+\bar u(t)$, and
\begin{equation}
l(x):=\rho(|x|)+\bar\varphi(0)+1,\ \ \RE\ x\in \R.
\end{equation}

{\bf Proof.}\ It follows from (2A1) and (3A2) that $\as$, for each $(y,z)\in \R\times \R^d$,
$$g(\omega,t,y,z)\ {\rm sgn}(y)\leq u(t)\rho(|y|)+|g(t,0,z)|\leq [u(t)+\bar u(t)]l(y)+\bar h(0)|z|^2,$$
where $l(x)$ is defined in (43). Note that $\rho$ is of linear growth. We can deduce that $l(\cdot)$ belongs to $\LS$ and then (3A1) holds true for the generator $g$. Thus, the existence of a bounded solution of BSDE$(\xi,g)$ follows from Theorem 3.1. Finally, the uniqueness follows directly from Theorems 2.1, 4.1 and 4.2. The proof is then completed.\vspace{0.2cm}\hfill $\Box$

{\bf Remark 4.1}\ By the proof of Theorem 2.1 and Theorems 4.1-4.3 we can conclude that if assumptions (2A2), (4A1) and (4A2) hold true only when $y\in [L_0,U_0]$, then the conclusions of Theorem 4.3 hold still true.

\section{Existence and uniqueness of $L^p\ (p>1)$ solutions}

In this section, we will fix a real number $p>1$ and establish two existence results, two comparison theorems and an existence and uniqueness result for the (maximal and minimal) $L^p$ solutions of BSDEs. As before, let us first introduce the following assumptions on the generator $g$, where we also assume that $0<T\leq +\infty$.\vspace{0.2cm}

{\bf (5A1)} There exist two functions $u(\cdot)\in \US$, $v(\cdot)\in \VS$ and a process $f_t\in \FS$ such that $\as$,
$$g(\omega,t,y,z)\ {\rm sgn}(y)\leq f_t(\omega)+u(t)|y|+v(t)|z|,\ \ \RE\ y,z.$$

{\bf (5A2)} There exist two functions $\bar u(\cdot)\in \US$, $\bar v(\cdot)\in \VS$ and a process $\bar f_t\in \FS$ such that $\as$,
$$\RE\ y\leq 0,\ \ \RE\ z\in \R^d,\ \ g(\omega,t,y,z)\leq \bar f_t(\omega)+\bar u(t)|y|+\bar v(t)|z|.$$

{\bf (5A3)} There exist two functions $\bar u(\cdot)\in \US$, $\bar v(\cdot)\in \VS$ and a process $\bar f_t\in \FS$ such that $\as$,
$$\RE\ y\geq 0,\ \ \RE\ z\in \R^d,\ \ -g(\omega,t,y,z)\leq \bar f_t(\omega)+\bar u(t)|y|+\bar v(t)|z|.$$

{\bf (5A4)} There exist two functions $u(\cdot)\in \US$, $v(\cdot)\in \VS$ and a process $f_t\in \FS$ such that $\as$,
$$|g(\omega,t,y,z)|\leq f_t(\omega)+u(t)|y|+v(t)|z|,\ \ \RE\ y,z.$$

{\bf (5A5)} There exists a nonnegative $(\F_t)$-progressively measurable process $(\check u_t)_{t\in\T}$ with $\ps,\ \int_0^T \check u_t(\omega)\ {\rm d}t<+\infty$, and two continuous functions $\check\varphi(\cdot),\ \check h(\cdot):\R\mapsto \R_+$ such that $\as$,
$$|g(\omega,t,y,z)|\leq \check u_t(\omega)\check \varphi(y)+\check h(y)|z|^2,\ \ \RE\ y,z.$$

{\bf Remark 5.1}\ It is clear that (5A4)$\Longleftrightarrow$(5A1)+(5A2)+(5A3), (5A4)$\Longrightarrow$(5A5), and (3A2)$\Longrightarrow$(5A5).\vspace{0.2cm}

The following lemma will be used in this section, which comes from Proposition 2.3 in \citet{FJ12}.\vspace{0.1cm}

{\bf Lemma 5.1}\ Assume that $0<T\leq +\infty$, $g$ is a generator of BSDEs and the process $g(t,0,0)\in \FS$. Assume further that there exist two functions $u(\cdot)\in \US$ and $v(\cdot)\in \VS$ such that $\as$,
\begin{equation}
|g(\omega,t,y_1,z_1)-g(\omega,t,y_2,z_2)|\leq u(t)|y_1-y_2|+ v(t)|z_1-z_2|,\ \ \RE\ y_1,y_2,z_1,z_2.
\end{equation}
Then for each $\xi\in \Lp$, BSDE$(\xi,g)$ has a unique $L^p$ solution.\vspace{0.2cm}

By Theorem 2.1 and Remark 2.2, we can establish the following lemma.\vspace{0.1cm}

{\bf Lemma 5.2} Assume that $0<T\leq +\infty$ and $g$ satisfies (2A3) and (5A1). Let $(y_t,z_t)_{t\in \T}$ be any $L^p$ solution of BSDE$(\xi,g)$ and $(y'_t,z'_t)_{t\in \T}$ the unique $L^p$ solution of BSDE$(|\xi|,g')$ by Lemma 5.1, where
$$g'(\omega,t,y,z):=f_t(\omega)+u(t)|y|+ v(t)|z|,\ \ \RE\ \omega,t,y,z.$$
Then for each $t\in \T$, we have
$$|y_t|\leq y'_t\ \ \ps.$$

{\bf Proof.}\ Note that $|\xi|\geq 0\ \ \ps$ and $g'(t,0,0)=f_t\geq 0\ \ \as$. By Theorem 2.1 we know that for each $t\in \T$,
\begin{equation}
y'_t\geq 0\ \ \ \ps.
\end{equation}
It follows from (5A1) that $\as$,
\begin{equation}
\RE\ y> 0,\ \RE\ z\in \R^d,\ \ g(\omega,t,y,z)\leq g'(\omega,t,y,z).
\end{equation}
Thus, in view of (45), (46), Remark 2.2 and the fact that $\xi\leq |\xi|\ \ \ps$, by Theorem 2.1 we deduce that for each $t\in \T$,
$$y_t\leq y'_t\ \ \ \ps.$$
Furthermore, by (5A1) we can also deduce that $\as$,
\begin{equation}
\RE\ y< 0,\ \RE\ z\in \R^d,\ \ -g'(\omega,t,y,z)\leq g(\omega,t,y,z).
\end{equation}
On the other hand, it is not difficult to verify that $(-y'_t,-z'_t)_{t\in \T}$ is the unique $L^p$ solution of BSDE$(-|\xi|,-g')$. Then, in view of (45), (47), Remark 2.2 and the fact that $-|\xi|\leq \xi\ \ \ps$, by Theorem 2.1 we can deduce that for each $t\in \T$,
$$-y'_t\leq y_t\ \ \ \ps.$$
Thus, we have completed the proof of Lemma 5.2. \vspace{0.2cm}\hfill $\Box$

The following Theorem 5.1 establishes an existence result for $L^p$ solutions of BSDEs, which is one of main results in this section. It improves virtually the corresponding existence results obtained respectively in \citet{Lep97}, \citet{Chen10} and \citet{Bri07} even for the case of the finite time horizon.\vspace{0.1cm}

{\bf Theorem 5.1} Assume that $0<T\leq +\infty$ and $g$ satisfies (2A3), (3A2) and (5A1) with $f_t\in \US$. Then for each $\xi\in \Lp$, BSDE$(\xi,g)$ has an $L^p$ solution.\vspace{0.1cm}

{\bf Proof.}\ Let us first assume that $\xi$ is nonnegative. Note first that (5A1) with $f_t\in \US$ can imply (3A1). Indeed, if $g$ satisfies (5A1), then $\as$,
$$g(\omega,t,y,z)\ {\rm sgn}(y)\leq \tilde u(t)\tilde l(y)+|z|^2,\ \ \RE\ y,z,\vspace{0.1cm}$$
where $\tilde u(t):=f_t+u(t)+v^2(t)\in\US$ and $\tilde l(x):=1+|x|\in \LS$.

It follows from Theorem 3.1 that BSDE$(\xi_n,g)$ has a maximal bounded solution $(y^n_t,z^n_t)_{t\in \T}$ for each $n\geq 1$, where $\xi_n:=\xi\wedge n$. Furthermore, in view of $\xi_n\leq \xi_{n+1}$, by Corollary 3.1 we also know that for each $t\in \T$, $(y^n_t)_{n=1}^{+\infty}$ is nondecreasing.

In the sequel, for each $(\omega,t,y,z)\in \Omega\times \T\times \R\times\R^d$, we define
$$g'(\omega,t,y,z):= f_t+u(t)|y|+v(t)|z|.$$
It follows from Lemma 5.1 that for each $n\geq 1$, BSDE$(|\xi_n|,g')$ has a unique $L^p$ solution $({}^ny'_t,{}^nz'_t)_{t\in \T}$. Furthermore, in view of (5A1), by Lemma 5.2 we can conclude that for each $t\in \T$ and each $n\geq 1$, $|y^n_t|\leq {}^ny'_t\ \ \ps$.

On the other hand, in view of $\xi\in \Lp$, it follows from Lemma 5.1 that BSDE$(|\xi|,g')$ has also a unique $L^p$ solution $(y'_t,z'_t)_{t\in \T}$. Then, in view of $|\xi_n|\leq |\xi|$, by Theorem 2.1 we know that for each $t\in \T$ and each $n\geq 1$, ${}^ny'_t\leq y'_t\ \ps$.

As a result, we have proved that for each $t\in \T$ and $n\geq 1$,
\begin{equation}
-y'_t\leq y^n_t\leq y^{n+1}_t\leq y'_t\ \ \ \ps.
\end{equation}
We define $y_\cdot=\Lim y^n_\cdot$, then
$$
\RE\ t\in[0,T],\ \ |y_t|\leq |y'_t|=y'_t \ \ \ \ps.\vspace{-0.2cm}
$$

In the sequel, we will use the localization procedure used in \citet{Bri06} to construct the desired solution. For each $k\geq 1$, let us introduce the following stopping time:
$$\tau_k=\inf\{t\in[0,T]: |y'_t|\geq k\}\wedge T.\vspace{0.1cm}$$
Then $(y^n_{k}(t),z^n_{k}(t)):=(y^n_{t\wedge \tau_k },z^n_t \mathbbm{1}_{ t\leq \tau_k})$ solves the
following BSDE:
\begin{equation}
  y_{k}^n(t)=y^n_{\tau_k}+\int_t^T \mathbbm{1}_{ s\leq \tau_k}
  g(s,y_{k}^n(s),z_{k}^n(s))
  {\rm d}s-\int_t^T z_{k}^n(s)\cdot {\rm d}B_s.
\end{equation}

It is very important to observe that $y_{k}^n$ is nondecreasing in $n$ and that, from the definition of $\tau_k$ and inequality (48),
$$
\sup\limits_{n\geq 1}\sup\limits_{t\in[0,T]}
\|y_{k}^n(t)\|_\infty\leq k.
$$
Furthermore, by (3A2) we know that $\as$, for each $(y,z)\in [-k,k]\times\R^d$,
$$|g(\omega,t,y,z)|\leq \bar u(t)\left(\max\limits_{x\in [-k,k]}\bar\varphi(x)\right)+\left(\max\limits_{x\in [-k,k]}\bar h(x)\right)|z|^2.$$
Thus, arguing as in the last second paragraph of the proof of Lemma 3.3, we can take the limit with respect to $n$ ($k$ being fixed) in (49) in the space $\s^2\times {\rm M}^2$. In particular, setting $y_k(t)=\sup_{n\geq 1}y^n_k(t)$, we know that $y_k(\cdot)$ is continuous and that there exists a process $z_{k}(t)\in {\M}^2$ such that $\Lim z^n_{k}(t)=z_{k}(t)$ in ${\M}^2$ and $(y_{k}(t),z_{k}(t))$ solves the following BSDE
\begin{equation}
  y_{k}(t)=\sup_{n\geq 1} y^n_{\tau_k}+\int_t^T \mathbbm{1}_{ s\leq \tau_k}
  g(s,y_{k}(s),z_{k}(s))
  {\rm d}s-\int_t^T z_{k}(s)\cdot {\rm d}B_s.
\end{equation}

Since $\tau_k\leq \tau_{k+1}$, it follows from the definitions of $y_{k}(\cdot),z_{k}(\cdot)$ and $y_\cdot$ that
$$y_{t\wedge\tau_k}=y_{k+1}(t\wedge\tau_k)=y_{k}(t)=\sup\limits_{n\geq 1}y^n_{t\wedge\tau_k},\ \ z_{k+1}(t)\mathbbm{1}_{ t\leq \tau_k}=z_{k}(t)=\lim\limits_{n\To \infty} z^n_t\mathbbm{1}_{ t\leq \tau_k}.$$
Thus, since $y_{k}(\cdot)$ are continuous processes and moreover $\ps,\  \tau_k=T$ for $k$ large enough, we know that $y_\cdot$ is continuous on $[0,T]$. Then we define $z_\cdot$ on $(0,T)$ by setting
$$z_t=z_{k}(t),\ \ {\rm if}\ t\in(0,\tau_k),$$
so that $z_t\mathbbm{1}_{ t\leq \tau_k}=z_k(t)\mathbbm{1}_{ t\leq \tau_k}=z_k(t)$ and (49) can be rewritten as
\begin{equation}
  y_{t\wedge\tau_k}=y_{\tau_k}
  +\int_{t\wedge\tau_k}^{\tau_k}
  g(s,y_s,z_s) {\rm d}s-\int_{t\wedge\tau_k}^{\tau_k}z_s\cdot {\rm d}B_s.
\end{equation}

Furthermore, we have
$$
\begin{array}{lll}
\Dis P \left(\int_0^T|z_s|^2{\rm d}s=\infty\right)& =& \Dis P \left(\int_0^T|z_s|^2{\rm d}s=\infty,\tau_k=T\right)\\
&&\Dis + P \left(\int_0^T|z_s|^2{\rm
d}s=\infty,\tau_k<T\right)\\
&\leq & \Dis P \left(\int_0^{\tau_k}|z_{k}(s)|^2{\rm d}s=\infty\right)+
P\left(\tau_k<T\right),
\end{array}$$
and we deduce, since $\tau_k\uparrow T$, that
$$\int_0^T |z_s|^2{\rm d}s<\infty\ \ \ \ps.$$
Thus, letting $k\To \infty$ in (51), we can deduce that $(y_\cdot,z_\cdot)$ is a solution of BSDE$(\xi,g)$. Furthermore, by (48) we know that $y_\cdot\in \s^p$, and then arguing as in Lemma 3.1 of \citet{Bri03}, in view of (5A1), we can deduce that $z_\cdot\in {\rm M}^p$. That is to say, $(y_\cdot,z_\cdot)$ is also an $L^p$ solution of BSDE$(\xi,g)$.

In the general case, we can use a double approximation: $\xi_{n,p}:=\xi^+\wedge n-\xi^-\wedge p$ as in \citet{Bri06} and \citet{Bri08}. The proof is then complete.\hfill \vspace{0.2cm}$\Box$

The following Theorem 5.2 gives a new existence result on the maximal and minimal $L^p$ solution of BSDEs. In view of Remark 5.1, we know that it generalizes Theorem 1 in \citet{FJT11} (see Lemma 2.2) even for the case of $L^2$ solutions.\vspace{0.1cm}

{\bf Theorem 5.2} Assume that $0<T\leq +\infty$ and $g$ satisfies (2A3), (5A1) and (5A5). Assume further that $g$ satisfies (5A2) (resp. (5A3)). Then for each $\xi\in \Lp$, BSDE$(\xi,g)$ has a maximal (resp. minimal) $L^p$ solution.\vspace{0.1cm}

{\bf Proof.} We only prove the case of the maximal solution, another case is similar. Assume now that $\xi\in \Lp$ and that $g$ satisfies (2A3), (5A1), (5A5) and (5A2). By (5A1) and (5A2) we know that $\as$,
$$g(\omega,t,y,z)\leq (f_t(\omega)+\bar f_t(\omega))+(u(t)+\bar u(t))|y|+(v(t)+\bar v(t))|z|,\ \ \RE\ y,z.$$
Then for each $n\geq 1$ and $(\omega,t,y,z)\in \Omega\times\T\times\R\times \R^d$, we can define
\begin{equation}
g_n(\omega,t,y,z):=\sup\limits_{(u,v)\in \R^{1+d}}
\{g(\omega,t,u,v)-n(u(t)+\bar u(t))|y-u|-n(v(t)+\bar v(t))|z-v|\}.
\end{equation}
Arguing as in \citet{FJT11}, we can also conclude that the sequence of functions $g_n$ is well defined for each $n\geq 1$, and it satisfies, $\as$,\vspace{0.2cm}

\ \ (i)\ $\RE\ y,z,\ g_n(\omega,t,y,z)\leq (f_t(\omega)+\bar f_t(\omega))+(u(t)+\bar u(t))|y|+(v(t)+\bar v(t))|z|.$

\ (ii)\ $\RE\ y,z,\ g_n(\omega,t,y,z)$ non-increases in $n$.

(iii)\ $\RE\ y_1,y_2,z_1,z_2$, we have
$$\hspace*{1.4cm}|g_n(\omega,t,y_1,z_1)-g_n(\omega,t,y_2,z_2)|\leq
n(u(t)+\bar u(t))|y_1-y_2|+n(v(t)+\bar v(t))|z_1-z_2|.\vspace{-0.1cm}$$

(iv)\ If $(y_n,z_n)\To (y,z)$, then $g_n(\omega,t,y_n,z_n)\to g(\omega,t,y,z)$.\vspace{0.2cm}

\noindent Furthermore, it follows from (52) and (5A1) that for each $n\geq 1$, $\as$,
\begin{equation}
\RE\ y<0,\ \ \RE\ z\in \R^d,\ \ g'(\omega,t,y,z)\leq g(\omega,t,y,z)\leq g_n(\omega,t,y,z),
\end{equation}
where
$$g'(\omega,t,y,z):=-f_t(\omega)-u(t)|y|-v(t)|z|,\ \ \RE\ \omega,t,y,z.$$

Note that (iii) and $|g_n(t,0,0)|\in \FS$ by (i) and (53). It follows from Lemma 5.1 that BSDE$(\xi,g_n)$ has a unique $L^p$ solution $(y^n_t,z^n_t)_{t\in \T}$ for each $n\geq 1$. In view of (ii) and (iii), by Theorem 2.1 we know that for each $t\in \T$, the sequel $(y^n_t)_{n=1}^{+\infty}$ is non-increasing. On the other hand, let $(y'_\cdot,z'_\cdot)$ be the unique $L^p$ solution of BSDE$(-|\xi|,g')$ by Lemma 5.1. Note that $y'_\cdot\leq 0$ by Theorem 2.1. In view of (53), Remark 2.2 and the fact that $-|\xi|\leq \xi\ \ \ps$, by Theorem 2.1 we know that for each $t\in \T$ and each $n\geq 1$,
\begin{equation}
y'_t\leq y^{n+1}_t\leq y^n_t\leq y^1_t\ \ \ \ps.\vspace{-0.1cm}
\end{equation}

In the sequel, for each $k\geq 1$, let ue introduce the following stopping time:
$$\tau_k=\inf\{t\in[0,T]: |y'_t|+|y^1_t|+\int_0^t (\check u_t+f_t+\bar f_t)\ {\rm d}s\geq k\}\wedge T.\vspace{0.1cm}$$
Then $(y^n_{k}(t),z^n_{k}(t)):=(y^n_{t\wedge \tau_k },z^n_t \mathbbm{1}_{ t\leq \tau_k})$ solves the
following BSDE:
$$
  y_{k}^n(t)=y^n_{\tau_k}+\int_t^T \mathbbm{1}_{ s\leq \tau_k}
  g_n(s,y_{k}^n(s),z_{k}^n(s))
  {\rm d}s-\int_t^T z_{k}^n(s)\cdot {\rm d}B_s.
$$
In view of (i) and the facts that $g_n\geq g$ and $g$ satisfies (5A5), we know that $\as$, for each $(y,z)\in [-k,k]\times\R^d$,
$$
\begin{array}{lll}
\Dis\sup\limits_{n\geq 1}|g_n(\omega,t,y,z)|&\leq &\Dis \check u_t(\omega)\left(\max\limits_{x\in [-k,k]}\check\varphi(x)\right)+f_t(\omega)+\bar f_t(\omega)+k(u(t)+\bar u(t))\\
&&\Dis+(v(t)+\bar v(t))^2+\left(\max\limits_{x\in [-k,k]}\check h(x)+1\right)|z|^2.
\end{array}
$$
Thus, in view of (54) and the definition of $\tau_k$, by a similar argument to Lemma 4.1 we can deduce that $(y^n_k(t),z^n_k(t))\in \s^{\infty}\times\M^2$ for each $n\geq 1$ and $k\geq 1$. Furthermore, in view of (iv), (54) and the above inequality together with the definition of $\tau_k$, arguing as in the proof of Theorem 5.1, we can define $y_\cdot=\lim_{n\To \infty} y^n_\cdot$ and use the localization procedure to obtain a process $z_\cdot$ such that $(y_\cdot,z_\cdot)$ is an $L^p$ solution of BSDE$(\xi,g)$.

Finally, it remains to show that $(y_\cdot,z_\cdot)$ is the maximal one among all $L^p$ solutions of BSDE$(\xi,g)$. Indeed, let $(\hat y_\cdot,\hat z_\cdot)$ be any $L^p$ solutions of BSDE$(\xi,g)$. In view of (iii) and the fact that $g_n\geq g$, by Theorem 2.1 we can obtain that for each $n\geq 1$ and each $t\in \T$,
$$\hat y_t\leq y^n_t\ \ \ \ps.$$
Letting $n\To \infty$ yields the desired result. The proof is then completed. \vspace{0.2cm} \hfill $\Box$

By Theorem 5.2 and Remark 5.1, the following corollary is immediate. It generalizes Lemma 2.2 to the case of $L^p$ solutions.\vspace{0.1cm}

{\bf Corollary 5.1}\ Assume that $0<T\leq +\infty$ and $g$ satisfies (2A3) and (5A4). Then for each $\xi\in \Lp$, BSDE$(\xi,g)$ has both a maximal $L^p$ solution and a minimal $L^p$ solution.

Similar to Theorems 2.2 and 2.3, by virtue of Theorem 2.1 we can prove the following Theorems 5.3 and 5.4, which establish the comparison theorems on the maximal and minimal $L^p$ solutions of BSDEs under the assumptions of Theorem 5.2. In view of Remark 5.1, they generalizes Theorems 2.2 and 2.3 even for the case of $L^2$ solutions.\vspace{0.1cm}

{\bf Theorem 5.3}\ \ Assume that $0<T\leq +\infty$, $g$ and $g'$ are two generators of BSDEs, and $(y_\cdot,z_\cdot)$ is any $L^p$ solution of BSDE$(\xi,g)$. Assume further that $g'$ satisfies (2A3), (5A1), (5A2) and (5A5), and $(y'_\cdot,z'_\cdot)$ is the maximal $L^p$ solution of BSDE$(\xi',g')$ by Theorem 5.2. If $\xi\leq \xi'\ \ps$ and (9) holds true, then for each $t\in\T$,
$$y_t\leq y'_t\ \ \ \ps.$$

{\bf Theorem 5.4}\ \ Assume that $0<T\leq +\infty$, $g$ and $g'$ are two generators of BSDEs, and $(y'_\cdot,z'_\cdot)$ is any $L^p$ solution of BSDE$(\xi',g')$. Assume further that $g$ satisfies (2A3), (5A1), (5A3) and (5A5), and $(y_\cdot,z_\cdot)$ is the minimal $L^p$ solution of BSDE$(\xi,g)$ by Theorem 5.2. If $\xi\leq \xi'\ \ps$ and (8) holds true, then for each $t\in\T$,
$$y_t\leq y'_t\ \ \ \ps.$$

{\bf Corollary 5.2} Assume that $0<T\leq +\infty$ and both $g$ and $g'$ satisfy (2A3) and (5A4). Let $(y_\cdot,z_\cdot)$ and $(y'_\cdot,z'_\cdot)$ be, respectively, the maximal (resp. minimal) $L^p$ solution of BSDE$(\xi,g)$ and BSDE$(\xi',g')$ by Corollary 5.1. If $\xi\leq \xi'\ \ps$ and for each $(y,z)\in \R\times \R^d$,
$$g(t,y,z)\leq g'(t,y,z)\ \ \ \as,\vspace{0.1cm}$$
then for each $t\in \T$, we have
$$y_t\leq y'_t\ \ \ \ps.\vspace{0.1cm}$$

Finally, by virtue of Theorem 5.1 and Theorem 2.1 we can obtain the following existence and uniqueness result for $L^p$ solutions of BSDEs.\vspace{0.1cm}

{\bf Theorem 5.5} Assume that $0<T\leq +\infty$ and that $g$ satisfies (2A1)-(2A3) and (3A2). Then for each $\xi\in \Lp$, BSDE$(\xi,g)$ has a unique $L^p$ solution.

{\bf Proof.}\ It follows from (2A1), (2A2) and (3A2) that $\as$, for each $(y,z)\in \R\times \R^d$,
\begin{equation}
\begin{array}{lll}
g(\omega,t,y,z)\ {\rm sgn}(y)&\leq & u(t)\rho(|y|)+|g(t,0,z)-g(t,0,0)|+|g(t,0,0)|\\
&\leq & u(t)(k|y|+k)+v(t)(a|z|+b)+\bar u(t)\bar\varphi(0)\\
&=& ku(t)+bv(t)+\bar u(t)\bar\varphi(0)+ku(t)|y|+av(t)|z|,
\end{array}
\end{equation}
where $k>0$ is the constant of linear growth for the function $\rho$. Consequently, (5A1) with $f_t=ku(t)+bv(t)+\bar u(t)\bar\varphi(0)\in \US$ holds true for the generator $g$. Thus, the existence of an $L^p$ solution of BSDE$(\xi,g)$ follows from Theorem 5.1. Finally, the uniqueness follows directly from Theorem 2.1. The proof is then completed.\hfill $\Box$

\section{Existence and uniqueness of $L^1$ solutions}

In this section, we will establish two existence results, two comparison theorems and an existence and uniqueness result for the (maximal and minimal) $L^1$ solutions of BSDEs. For convenience of expression, let us first fix a real number $\alpha\in (0,1)$ and let $\vS$ denote the set of nonnegative functions $\lambda(\cdot)$ from $\T$ to $\R_+$ such that
$$\int_0^T \lambda^{2\over 2-\alpha}(t)\ {\rm d}t<+\infty.\vspace{0.2cm}$$
In the sequel, let us introduce the following assumptions on the generator $g$, where $0<T\leq +\infty$, and $a\wedge b$ represents the minimal number between $a$ and $b$.\vspace{0.2cm}

{\bf (6A1)} There exists a nonnegative $(\F_t)$-progressively measurable process $(f_t)_{t\in\T}\in L^1(\T\tim\Omega)$ and three functions $u(\cdot)\in \US$, $v(\cdot)\in \VS$ and $\lambda(\cdot)\in \vS$ such that $\as$,
$$g(\omega,t,y,z)\ {\rm sgn}(y)\leq f_t(\omega)+u(t)|y|+(v(t)|z|)\wedge(\lambda(t)|z|^\alpha),\ \ \RE\ y,z.$$

{\bf (6A2)} There exists a nonnegative $(\F_t)$-progressively measurable process $(\bar f_t)_{t\in\T}\in L^1(\T\tim\Omega)$ and three functions $\bar u(\cdot)\in \US$, $\bar v(\cdot)\in \VS$ and $\bar\lambda(\cdot)\in \vS$ such that $\as$,
$$\RE\ y\leq 0,\ \ \RE\ z\in \R^d,\ \ g(\omega,t,y,z)\leq \bar f_t(\omega)+\bar u(t)|y|+(\bar v(t)|z|)\wedge(\bar\lambda(t)|z|^\alpha).$$

{\bf (6A3)} There exists a nonnegative $(\F_t)$-progressively measurable process $(\bar f_t)_{t\in\T}\in L^1(\T\tim\Omega)$ and three functions $\bar u(\cdot)\in \US$, $\bar v(\cdot)\in \VS$ and $\bar\lambda(\cdot)\in \vS$ such that $\as$,
$$\RE\ y\geq 0,\ \ \RE\ z\in \R^d,\ \ -g(\omega,t,y,z)\leq \bar f_t(\omega)+\bar u(t)|y|+(\bar v(t)|z|)\wedge(\bar\lambda(t)|z|^\alpha).$$

{\bf (6A4)} There exists a nonnegative $(\F_t)$-progressively measurable process $(f_t)_{t\in\T}\in L^1(\T\tim\Omega)$ and three functions $u(\cdot)\in \US$, $v(\cdot)\in \VS$ and $\lambda(\cdot)\in \vS$ such that $\as$,
$$|g(\omega,t,y,z)|\leq f_t(\omega)+u(t)|y|+(v(t)|z|)\wedge(\lambda(t)|z|^\alpha),\ \ \RE\ y,z.$$

{\bf Remark 6.1}\ It is clear that (6A4)$\Longleftrightarrow$ (6A1)+(6A2)+(6A3), and (6A4)$\Longrightarrow$(5A5).\vspace{0.2cm}

The following lemma will be our basic tool in the treatment of $L^1$-solutions.\vspace{0.1cm}

{\bf Lemma 6.1}\ Assume that $0<T\leq +\infty$, the generator $g$ satisfies (44) and (2A5) (or (2A5')), and the process $g(t,0,0)\in L^1(\Omega\times\T)$. Then for each $\xi\in L^1(\Omega,\F_T,P)$, BSDE$(\xi,g)$ has a unique $L^1$ solution.\vspace{0.1cm}

{\bf Proof.} The uniqueness follows from Theorem 2.4. Using a similar argument to Theorem 6.3 in \citet{Bri03} and making use of Lemma 5.1 and (21), we can prove the existence. The proof is standard, we omit it here.\vspace{0.2cm} \hfill $\Box$

Similar to Lemma 5.2, we can prove the following Lemma 6.2.\vspace{0.1cm}

{\bf Lemma 6.2} Assume that $0<T\leq +\infty$ and the generator $g$ satisfies (2A3) and (6A1). Let $(y_t,z_t)_{t\in \T}$ be any $L^1$ solution of BSDE$(\xi,g)$ and $(y'_t,z'_t)_{t\in \T}$ the unique $L^1$ solution of BSDE$(|\xi|,g')$ by Lemma 6.1, where
$$g'(\omega,t,y,z):=f_t(\omega)+u(t)|y|+(v(t)|z|)
\wedge(\lambda(t)|z|^\alpha),\ \ \RE\ \omega,t,y,z.$$
Then for each $t\in \T$, we have\vspace{-0.1cm}
$$|y_t|\leq y'_t\ \ \ps.\vspace{-0.1cm}$$

{\bf Proof.}\ We consider the function $\kappa(t,x):=(v(t)x)\wedge(\lambda(t)x^\alpha)$ for $(t,x)\in \T\times\R_+$. It is not hard to verify that for each $t\in\T$, $\kappa(t,\cdot)$ is nondecreasing and sub-additive on $\R_+$, i.e., $\kappa(t,x_1+x_2)\leq \kappa(t,x_1)+\kappa(t,x_2)$ for each $x_1,x_2\in \R_+$. Based on this fact, we can prove that for each $t\in \T$ and $(z_1,z_2)\in \R^d\times\R^d$,
\begin{equation}
|(v(t)|z_1|)\wedge(\lambda(t)|z_1|^\alpha)-(v(t)|z_2|)\wedge(\lambda(t)
|z_2|^\alpha)|\leq (v(t)|z_1-z_2|)\wedge(\lambda(t)|z_1-z_2|^\alpha).
\end{equation}
Then we can see that $g'$ satisfies (44) and (2A5'). Thus, note that $|\xi|\geq 0\ \ \ps$ and $g'(t,0,0)=f_t\geq 0\ \ \as$, Theorem 2.4 yields that for each $t\in \T$,
\begin{equation}
y'_t\geq 0\ \ \ \ps.
\end{equation}
It follows from (6A1) that $\as$,
\begin{equation}
\RE\ y> 0,\ \RE\ z\in \R^d,\ \ g(\omega,t,y,z)\leq g'(\omega,t,y,z).
\end{equation}
Then, in view of (57), (58), Remark 2.2 and the fact that $\xi\leq |\xi|\ \ \ps$, by Theorem 2.4 we deduce that for each $t\in \T$,
$$y_t\leq y'_t\ \ \ \ps.$$
Furthermore, by (6A1) we can also deduce that $\as$,
\begin{equation}
\RE\ y< 0,\ \RE\ z\in \R^d,\ \ -g'(\omega,t,y,z)\leq g(\omega,t,y,z).
\end{equation}
On the other hand, it is not difficult to verify that $(-y'_t,-z'_t)_{t\in \T}$ is the unique $L^1$ solution of BSDE$(-|\xi|,-g')$. Then, in view of (57), (59), Remark 2.2 and the fact that $-|\xi|\leq \xi\ \ \ps$, by Theorem 2.4 we can deduce that for each $t\in \T$,
$$-y'_t\leq y_t\ \ \ \ps.$$
Thus, we have completed the proof of Lemma 6.2. \vspace{0.2cm}\hfill $\Box$

The following Theorem 6.1 establishes an existence result on $L^1$ solutions of BSDEs, which is one of main results in this section. It improves virtually the corresponding existence results obtained in \citet{Bri06} and \citet{Bri03} for the one dimensional case, even for the case of the finite time horizon.\vspace{0.1cm}

{\bf Theorem 6.1} Assume that $0<T\leq +\infty$ and $g$ satisfies (2A3), (3A2) and (6A1) with $f_t\in \US$. Then for each $\xi\in L^1(\Omega,\F_T,P)$, BSDE$(\xi,g)$ has an $L^1$ solution.\vspace{0.1cm}

{\bf Proof.}\ Let us first assume that $\xi$ is nonnegative. Note that (6A1)$\Longrightarrow$(5A1). By Theorem 3.1 and the proof of Theorem 5.1 we know that BSDE$(\xi_n,g)$ has a maximal bounded solution $(y^n_t,z^n_t)_{t\in \T}$ for each $n\geq 1$, where $\xi_n:=\xi\wedge n$. Furthermore, by Corollary 3.1 we also know that for each $t\in \T$, $(y^n_t)_{n=1}^{+\infty}$ is nondecreasing.

In the sequel, for each $(\omega,t,y,z)\in \Omega\times \T\times \R\times\R^d$, define
$$g'(\omega,t,y,z):= f_t+u(t)|y|+(v(t)|z|)\wedge (\lambda(t)|z|^\alpha).$$
By (56) we know that $g'$ satisfies (44) and (2A5'). It then follows from Lemma 6.1 that for each $n\geq 1$, BSDE$(|\xi_n|,g')$ has a unique $L^1$ solution $({}^ny'_t,{}^nz'_t)_{t\in \T}$. Furthermore, in view of (6A1), by Lemma 6.2 we can conclude that for each $t\in \T$ and each $n\geq 1$, $|y^n_t|\leq {}^ny'_t\ \ \ps$.

On the other hand, in view of $\xi\in L^1(\Omega,\F_T,P)$, it follows from Lemma 6.1 that BSDE$(|\xi|,g')$ has also a unique $L^1$ solution $(y'_t,z'_t)_{t\in \T}$. Then, in view of $|\xi_n|\leq |\xi|$, by Theorem 2.4 we know that for each $t\in \T$ and each $n\geq 1$, ${}^ny'_t\leq y'_t\ \ps$.

As a result, we have proved that for each $t\in \T$ and $n\geq 1$,
\begin{equation}
-y'_t\leq y^n_t\leq y^{n+1}_t\leq y'_t\ \ \ps.
\end{equation}
In the sequel, arguing as in Theorem 5.1, we define $y_\cdot=\lim_{n\To\infty} y^n_\cdot$ and use the localization procedure to find a process $z_\cdot$ such that $(y_\cdot,z_\cdot)$ is a solution of BSDE$(\xi,g)$. Furthermore, by (60) we know that $y_\cdot$ belongs to the class (D) and the space  $\s^\beta$ for each $\beta\in (0,1)$, and then arguing as in Lemma 3.1 of \citet{Bri03}, in view of (6A1), we can deduce that $z_\cdot\in {\rm M}^\beta$ for each $\beta\in (0,1)$. That is to say, $(y_\cdot,z_\cdot)$ is also an $L^1$ solution of BSDE$(\xi,g)$.

In the general case, we can use a double approximation: $\xi_{n,p}:=\xi^+\wedge n-\xi^-\wedge p$ as in \citet{Bri06} and \citet{Bri08}. The proof is then complete.\hfill \vspace{0.2cm}$\Box$

The following Theorem 6.2 gives a new existence result on the maximal and minimal $L^1$ solution of BSDEs. \vspace{0.1cm}

{\bf Theorem 6.2} Assume that $0<T\leq +\infty$ and $g$ satisfies (2A3), (6A1) and (5A5). Assume further that $g$ satisfies (6A2) (resp. (6A3)). Then for each $\xi\in L^1(\Omega,\F_T,P)$, BSDE$(\xi,g)$ has a maximal (resp. minimal) $L^1$ solution.\vspace{0.1cm}

{\bf Proof.} We only prove the case of the maximal solution, another case is similar. Assume now that $\xi\in L^1(\Omega,\F_T,P)$ and that $g$ satisfies (2A3), (6A1), (6A2) and (5A5). By (6A1) and (6A2) we know that $\as$,
$$g(\omega,t,y,z)\leq \tilde f_t(\omega)+\tilde u(t)|y|+(\tilde v(t)|z|)\wedge(\tilde \lambda(t)|z|^\alpha),\ \ \RE\ y,z.$$
where $\tilde f_t(\omega):=f_t(\omega)+\bar f_t(\omega)\in L^1(\T\tim\Omega)$, $\tilde u(t):=u(t)+\bar u(t)\in \US$, $\tilde v(t):=v(t)+\bar v(t)\in \VS$ and $\tilde \lambda(t):=\lambda(t)+\bar \lambda(t)\in \vS$. Then for each $n\geq 1$ and $(\omega,t,y,z)\in \Omega\times\T\times\R\times \R^d$, we can define
\begin{equation}
g_n(\omega,t,y,z):=\sup\limits_{(u,v)\in \R^{1+d}}
\{g(\omega,t,u,v)-n\tilde u(t)|y-u|-n[(\tilde v(t)|z-v|)\wedge(\tilde \lambda(t)|z-v|^\alpha)]\}.
\end{equation}
In view of (56), arguing as in \citet{FJT11}, we can also prove that the sequence of functions $g_n$ is well defined for each $n\geq 1$, and it satisfies, $\as$,\vspace{0.2cm}

\ \ (i)\ $\RE\ y,z,\ g_n(\omega,t,y,z)\leq \tilde f_t(\omega)+\tilde u(t)|y|+(\tilde v(t)|z|)\wedge(\tilde \lambda(t)|z|^\alpha).$

\ (ii)\ $\RE\ y,z,\ g_n(\omega,t,y,z)$ non-increases in $n$.

(iii)\ $\RE\ y_1,y_2,z_1,z_2$, we have
$$\hspace*{1.4cm}|g_n(\omega,t,y_1,z_1)-g_n(\omega,t,y_2,z_2)|\leq
n\tilde u(t)|y_1-y_2|+n[(\tilde v(t)|z_1-z_2|)\wedge(\tilde \lambda(t)|z_1-z_2|^\alpha)].\vspace{-0.1cm}$$

(iv)\ If $(y_n,z_n)\To (y,z)$, then $g_n(\omega,t,y_n,z_n)\to g(\omega,t,y,z)$.\vspace{0.2cm}

\noindent Furthermore, it follows from (61) and (6A1) that for each $n\geq 1$, $\as$,
\begin{equation}\RE\ y<0,\ \ \RE\ z\in \R^d,\ \ g'(\omega,t,y,z)\leq g(\omega,t,y,z)\leq g_n(\omega,t,y,z),
\end{equation}
where
$$g'(\omega,t,y,z):=-f_t(\omega)-u(t)|y|-(v(t)|z|)\wedge(\lambda(t)
|z|^\alpha),\ \ \RE\ \omega,t,y,z.\vspace{0.1cm}$$
It follows from (56) that $g'$ satisfies (44) and (2A5').

Note that (iii) can imply that $g_n$ satisfies (44) and (2A5'), and note that $|g_n(t,0,0)|\in L^1(\Omega\times\T)$ by (i) and (62). It follows from Lemma 6.1 that BSDE$(\xi,g_n)$ has a unique $L^1$ solution $(y^n_t,z^n_t)_{t\in \T}$ for each $n\geq 1$. In view of (ii) and (iii), by Theorem 2.4 we know that for each $t\in \T$, the sequel $(y^n_t)_{n=1}^{+\infty}$ is non-increasing. On the other hand, by Lemma 6.1 we can let $(y'_\cdot,z'_\cdot)$ be the unique $L^1$ solution of BSDE$(-|\xi|,g')$. Note that $y'_\cdot\leq 0$ by Theorem 2.4. In view of (62), (iii), Remark 2.2 and the fact that $-|\xi|\leq \xi\ \ \ps$, by Theorem 2.4 we know that for each $t\in \T$ and each $n\geq 1$,
\begin{equation}
y'_t\leq y^{n+1}_t\leq y^n_t\leq y^1_t\ \ \ \ps.
\end{equation}
In the sequel, arguing as in the proof of Theorems 5.2 and 6.1, in view of (63) and (iv), we can define $y_\cdot=\lim_{n\To \infty} y^n_\cdot$ and use the localization procedure to obtain a process $z_\cdot$ such that $(y_\cdot,z_\cdot)$ is an $L^1$ solution of BSDE$(\xi,g)$.

Finally, it remains to show that $(y_\cdot,z_\cdot)$ is the maximal one among all $L^1$ solutions of BSDE$(\xi,g)$. Indeed, let $(\hat y_\cdot,\hat z_\cdot)$ be any $L^1$ solutions of BSDE$(\xi,g)$. In view of (iii) and the fact that $g_n\geq g$, by Theorem 2.4 we can obtain that for each $n\geq 1$ and each $t\in \T$,
$$\hat y_t\leq y^n_t\ \ \ \ps.$$
Letting $n\To \infty$ yields the desired result. The proof is then completed. \vspace{0.2cm} \hfill $\Box$

By Theorem 6.2 and Remark 6.1, the following corollary is immediate.\vspace{0.1cm}

{\bf Corollary 6.1}\ Assume that $0<T\leq +\infty$ and the generator $g$ satisfies (2A3) and (6A4). Then for each $\xi\in L^1(\Omega,\F_T,P)$, BSDE$(\xi,g)$ has both a maximal $L^1$ solution and a minimal $L^1$ solution.\vspace{0.1cm}

{\bf Remark 6.2}\ A similar result to Corollary 6.1 has been obtained in the first version of \citet{Bri06}, where only is the case of $0<T<+\infty$ considered, all of $f_t$, $u(t)$, $v(t)$ and $\lambda(t)$ are constants and we do not know whether the $L^1$ solution constructed by them is the maximal (or minimal) one or not. Hence, Corollary 6.1 improves it.\vspace{0.2cm}

Similar to Theorems 2.2 and 2.3, by virtue of Theorem 2.4 we can prove the following Theorems 6.3 and 6.4, which establish the comparison theorems on the maximal and minimal $L^1$ solutions of BSDEs under the assumptions of Theorem 6.2.\vspace{0.1cm}

{\bf Theorem 6.3}\ \ Assume that $0<T\leq +\infty$, $g$ and $g'$ are two generators of BSDEs, and $(y_\cdot,z_\cdot)$ is any $L^1$ solution of BSDE$(\xi,g)$. Assume further that $g'$ satisfies (2A3), (6A1), (6A2) and (5A5), and $(y'_\cdot,z'_\cdot)$ is the maximal $L^1$ solution of BSDE$(\xi',g')$ by Theorem 6.2. If $\xi\leq \xi'\ \ps$ and (9) holds true, then for each $t\in\T$,
$$y_t\leq y'_t\ \ \ \ps.$$

{\bf Theorem 6.4}\ \ Assume that $0<T\leq +\infty$, $g$ and $g'$ are two generators of BSDEs, and $(y'_\cdot,z'_\cdot)$ is any $L^1$ solution of BSDE$(\xi',g')$. Assume further that $g$ satisfies (2A3), (6A1), (6A3) and (5A5), and $(y_\cdot,z_\cdot)$ is the minimal $L^1$ solution of BSDE$(\xi,g)$ by Theorem 6.2. If $\xi\leq \xi'\ \ps$ and (8) holds true, then for each $t\in\T$,
$$y_t\leq y'_t\ \ \ \ps.$$

{\bf Corollary 6.2} Assume that $0<T\leq +\infty$ and both $g$ and $g'$ satisfy (2A3) and (6A4). Let $(y_\cdot,z_\cdot)$ and $(y'_\cdot,z'_\cdot)$ be, respectively, the maximal (resp. minimal) $L^1$ solution of BSDE$(\xi,g)$ and BSDE$(\xi',g')$ by Corollary 6.1. If $\xi\leq \xi'\ \ps$ and for each $(y,z)\in \R\times \R^d$,
$$g(t,y,z)\leq g'(t,y,z)\ \ \ \as,\vspace{0.1cm}$$
then for each $t\in \T$, we have
$$y_t\leq y'_t\ \ \ \ps.$$

Finally, by virtue of Theorem 6.1 and Theorem 2.4 we can obtain the following new existence and uniqueness result for $L^1$ solutions.\vspace{0.1cm}

{\bf Theorem 6.5} Assume that $0<T\leq +\infty$ and that $g$ satisfies (2A1)-(2A3) and (3A2). Assume further that $g$ satisfies (2A5) with $f_t\in \US$ or (2A5'). Then for each $\xi\in L^1(\Omega,\F_T,P)$, BSDE$(\xi,g)$ has a unique $L^1$ solution.

{\bf Proof.}\ We only prove the case where $g$ satisfies (2A1)-(2A3), (3A2) and (2A5) with $f_t\in \US$. Another case is similar. It follows from (2A1), (2A5) with $f_t\in \US$, and (3A2) that $\as$, for each $(y,z)\in \R\times \R^d$,
\begin{equation}
\begin{array}{lll}
g(\omega,t,y,z)\ {\rm sgn}(y)&\leq & u(t)\rho(|y|)+|g(t,0,z)-g(t,0,0)|+|g(t,0,0)|\\
&\leq & u(t)(k|y|+k)+\lambda(t)(f_t+|z|)^\alpha+\bar u(t)\bar\varphi(0)\\
&\leq & \tilde f_t+ku(t)|y|+\lambda(t)|z|^\alpha,
\end{array}
\end{equation}
where $\tilde f_t:= ku(t)+\lambda(t)f_t^\alpha+\bar u(t)\bar\varphi(0)\in \US$ by H\"{o}lder inequality and the assumptions of $u(t), \lambda(t), f_t$ and $\bar u(t)$, and $k>0$ is the constant of linear growth for $\rho$. Thus, combining (55) and (64) yields that (6A1) holds true for the generator $g$. Thus, the existence of an $L^1$ solution of BSDE$(\xi,g)$ follows from Theorem 6.1. Finally, the uniqueness follows directly from Theorem 2.4. The proof is then completed.\hfill $\Box$\vspace{0.4cm}




\end{document}